\documentclass[12pt]{amsart}%{article}
\usepackage{amsmath}
\usepackage{amssymb}
\usepackage{amsthm}
\usepackage{graphicx}
\usepackage{mathabx}
\usepackage{mathrsfs}
\usepackage{enumitem}
\usepackage{comment}
\usepackage{soul}

\newtheorem{add}[equation]{Addendum}

\newtheorem{cor}[equation]{Corollary}

\newtheorem{lem}[equation]{Lemma}
\newtheorem{lemma}[equation]{Lemma}
\newtheorem{prop}[equation]{Proposition}

\newtheorem{thm}[equation]{Theorem}
\newtheorem{theorem}[equation]{Theorem}

\newtheorem{remark}[equation]{Remark}
\newtheorem{rem}[equation]{Remark}

\theoremstyle{bolddefinition}

\newtheorem{definition}[equation]{Definition}

\def\A{{\mathcal A}}
\def\B{{\mathrm B}}
\def\C{{\mathbb C}}
\def\R{{\mathbb R}}
\def\H{{\mathbb H}}

\def\N{{\mathbb N}}

\def\Z{{\mathbb Z}}
\def\Q{{\mathbb Q}}

\def\al{\alpha}
\def\ga{\gamma}
\def\Ga{\Gamma}

\def\eps{\epsilon}
\def\la{\lambda}
\def\La{\Lambda}
\def\si{\sigma}

\def\Om{\Omega}
\def\V{{\mathcal V}}

\def\geo{\partial_{\infty}}
\def\Isom{\mathop{\hbox{Isom}}}

\usepackage[pagebackref=true, colorlinks]{hyperref}
\hypersetup{pdffitwindow=true,linkcolor=blue,citecolor=blue,urlcolor=blue,menucolor=blue}

\usepackage{tikz}

\newcommand{\foh}{\tfrac12}

\newcommand{\dd}{\partial}

\newcommand{\bH}{\mathbb{H}}
\newcommand{\bS}{\mathbb{S}}
\newcommand{\bb}{\mathbf{b}}
\newcommand{\G}{\Gamma}      
\newcommand{\g}{\gamma}      
\def\cF{{\mathcal F}}
\def\cG{{\mathcal G}}

\def\al{\alpha}

\def\vep{\varepsilon}

\newcommand{\sP}{\mathscr{P}}
\newcommand{\sB}{\mathscr{B}}
\newcommand{\<}{\left\langle}
\renewcommand{\>}{\right\rangle}

\def\PSL{\mathop{\hbox{PSL}}}
\newcommand{\cO}{\mathcal{O}}
\newcommand{\fo}{\mathfrak{o}}
\newcommand{\cK}{\mathcal{K}}

\newtheorem{example}[equation]{Example}
\newtheorem{rmk}[equation]{Remark}

\newcommand\be{\begin{equation}}
\newcommand\ee{\end{equation}}
\newcommand\bp{\begin{pmatrix}}
\newcommand\ep{\end{pmatrix}}
\newcommand\pf{\begin{proof}}
\newcommand\epf{\end{proof}}

\newcommand{\appref}[1]{\hyperref[#1]{Appendix}}
\newcommand{\thmref}[1]{\hyperref[#1]{Theorem \ref{#1}}}
\newcommand{\propref}[1]{\hyperref[#1]{Proposition \ref{#1}}}
\newcommand{\corref}[1]{\hyperref[#1]{Corollary \ref{#1}}}
\newcommand{\defref}[1]{\hyperref[#1]{Definition \ref{#1}}}
\newcommand{\secref}[1]{\hyperref[#1]{\S\ref{#1}}}
\newcommand{\rmkref}[1]{\hyperref[#1]{Remark \ref{#1}}}
\newcommand{\figref}[1]{\hyperref[#1]{Figure \ref{#1}}}
\newcommand{\lemref}[1]{\hyperref[#1]{Lemma \ref{#1}}}
\newcommand{\exref}[1]{\hyperref[#1]{Example \ref{#1}}}
\newcommand{\addref}[1]{\hyperref[#1]{Addendum \ref{#1}}}

\newcommand{\mattwo}[4]
{\left(\begin{array}{cc}
                        #1  & #2   \\
                        #3 &  #4
                          \end{array}\right) }

\usepackage{tikz}
\usetikzlibrary{arrows,positioning} 
\tikzset{
    >=stealth',
    pil/.style={
           ->,
           thick,
           shorten <=2pt,
           shorten >=2pt,}
}

\begin{document}

\author{Michael Kapovich}
\thanks{Kapovich is supported by the NSF grant DMS-16-04241 and Simons Fellowship, grant number 391602.}
\email{kapovich@math.ucdavis.edu}
\address{Department of Mathematics, UC Davis, Davis, CA 95616}
\author{Alex Kontorovich}
\thanks{Kontorovich is supported by
an NSF CAREER grant DMS-1455705, an NSF FRG grant DMS-1463940,    NSF grant DMS-1802119, a BSF grant number 2014099, and the Simons Foundation through MoMath's Distinguished Visiting Professorship for the
Public Dissemination of Mathematics.}
\email{alex.kontorovich@rutgers.edu}
\address{Rutgers University, New Brunswick, NJ and National Museum of Mathematics, NY, NY}

\title[On Kleinian  Packings and Arithmetic Groups]{On SuperIntegral Kleinian Sphere Packings, Bugs, and Arithmetic Groups}
\date{\today}
\begin{abstract}
We develop the notion of a Kleinian Sphere Packing, a generalization of 
 ``crystallographic'' (Apollonian-like) sphere packings defined by Kontorovich-Nakamura \cite{KN}.
Unlike crystallographic packings, Kleinian packings exist in all dimensions, as do ``superintegral'' such. 
We extend the Arithmeticity Theorem to Kleinian packings, that is, the superintegral ones come from $\Q$-arithmetic lattices of simplest type.
The same holds for more general objects we call  Kleinian Bugs, in which the spheres need not be disjoint but can meet with dihedral angles $\pi/m$ for finitely many $m$. 
We settle two questions from \cite{KN}: $(i)$ that the Arithmeticity Theorem is in general false over number fields, and $(ii)$
 that  integral packings only arise from non-uniform lattices.
\end{abstract}

\maketitle

\section{Introduction}

The classical Apollonian packing in the plane, usually described by an ad hoc construction involving inscribing tangent circles, exhibits a number of thereafter surprising arithmetic and dynamical properties; see, e.g., \cite{Kontorovich2013}. 
In this paper, we complete the program initiated in \cite{KN} to understand the relationship between such  packings and the theory of 
%$\Q$-
arithmetic groups in hyperbolic space.

\subsection{Kleinian (and Crystallographic) Packings}\ % and Summary of \cite{KN}}\

\medskip
 A {\bf sphere packing} (or just ``packing'') $\sP$ of $\bS^n\cong\geo\bH^{n+1}$ ($n\ge2$) is an infinite collection of round balls in $\bS^n$ with pairwise disjoint interiors, such that the union of the balls is dense in $\bS^n$. 
 %({\color{red}\st{Throughout the paper, we will work with the upper half-space model of  hyperbolic space $\bH^{n+1}$. Accordingly, we will}} 
 We identify the ideal boundary $\bS^n=\dd_\infty\H^{n+1}$ with the one-point compactification  of Euclidean $n$-space $\bS^n=\R^n\cup\{\infty\}$.  By abuse of terminology, we will conflate the collection of balls $\sP$ with the collection of round spheres bounded by these balls. In view of the density condition, a packing contains balls with arbitrarily small Euclidean radii. (Of course, the radii are only defined once we choose an identification $\bS^n\cong \geo\bH^{n+1}$, in particular, we choose a point at $\infty$.) 
 The {\bf bend}\footnote{Note that in the theory of Kleinian groups, ``bend'' more often refers to dihedral angle; but for integral sphere packings, bend is used for inversive radii. For circles, bend is the curvature, but for higher dimensional spheres, (Gaussian) curvature is inverse square-radius.} of a sphere is the reciprocal of its radius, with the convention that  a sphere containing $\infty$ in its interior has negative radius and  bend.
A  packing is {\bf integral} if all its spheres have integer bends.

\begin{figure}
\includegraphics[width=.3\textwidth]{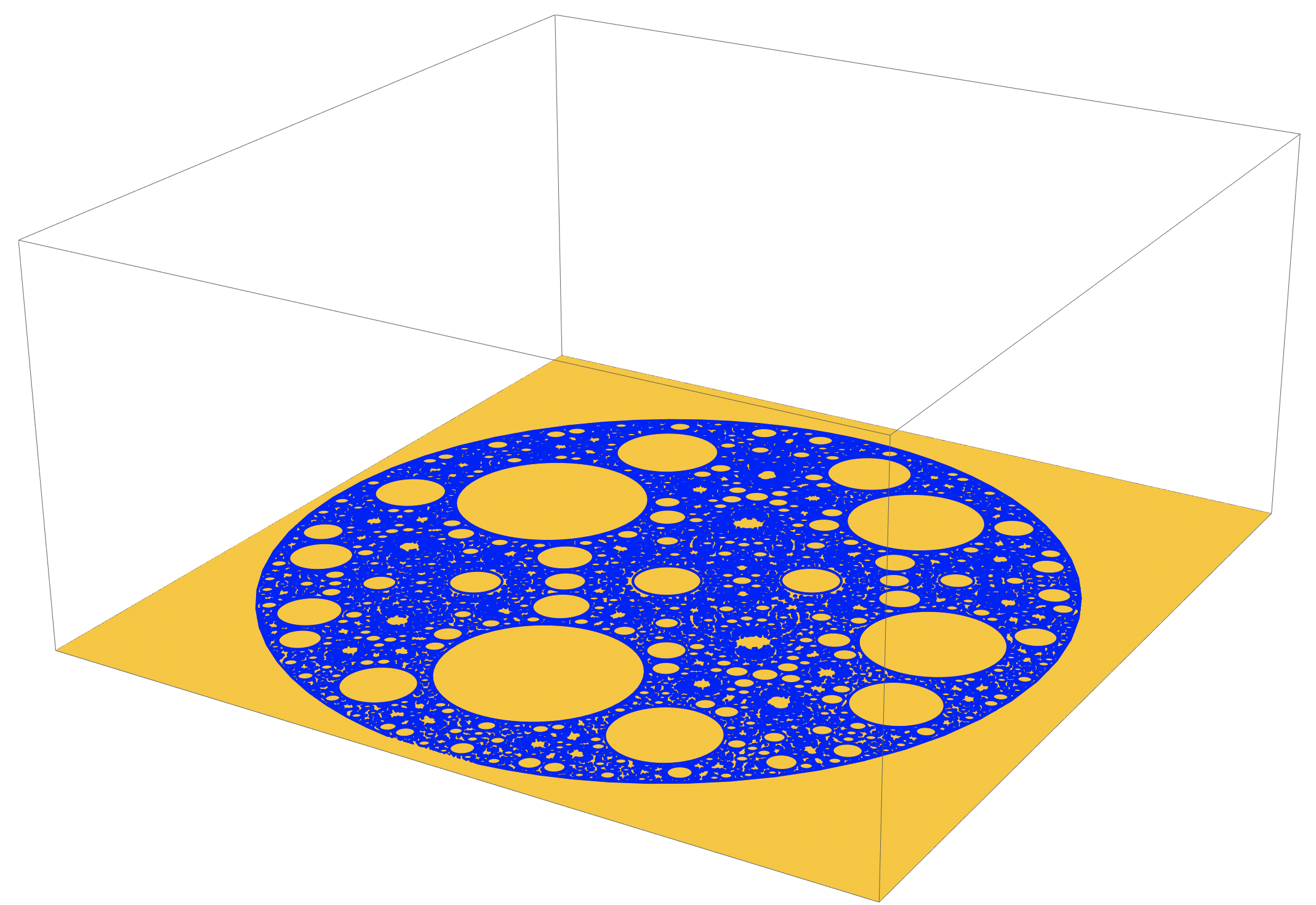}
\
\includegraphics[width=.3\textwidth]{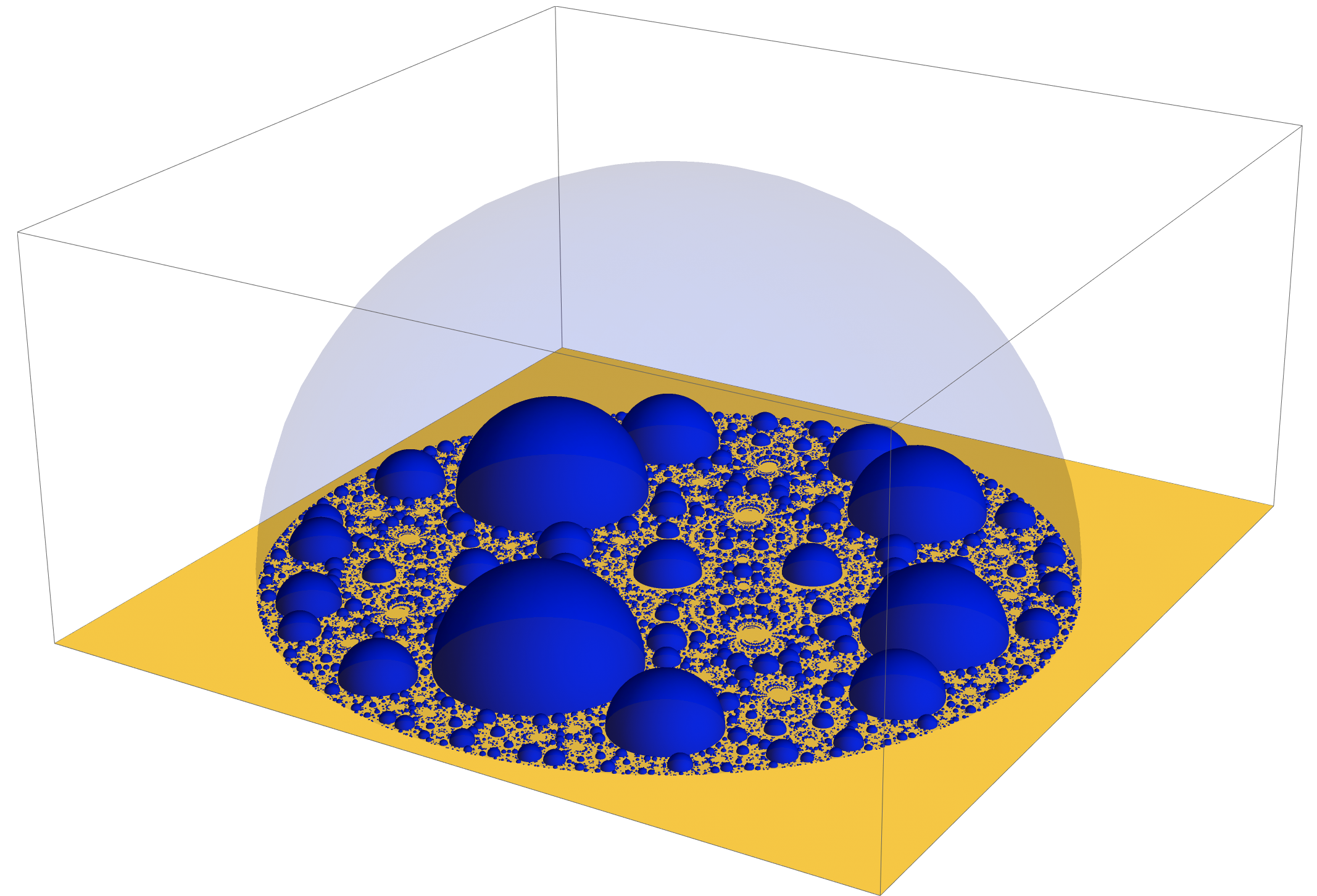}
\
\includegraphics[width=.3\textwidth]{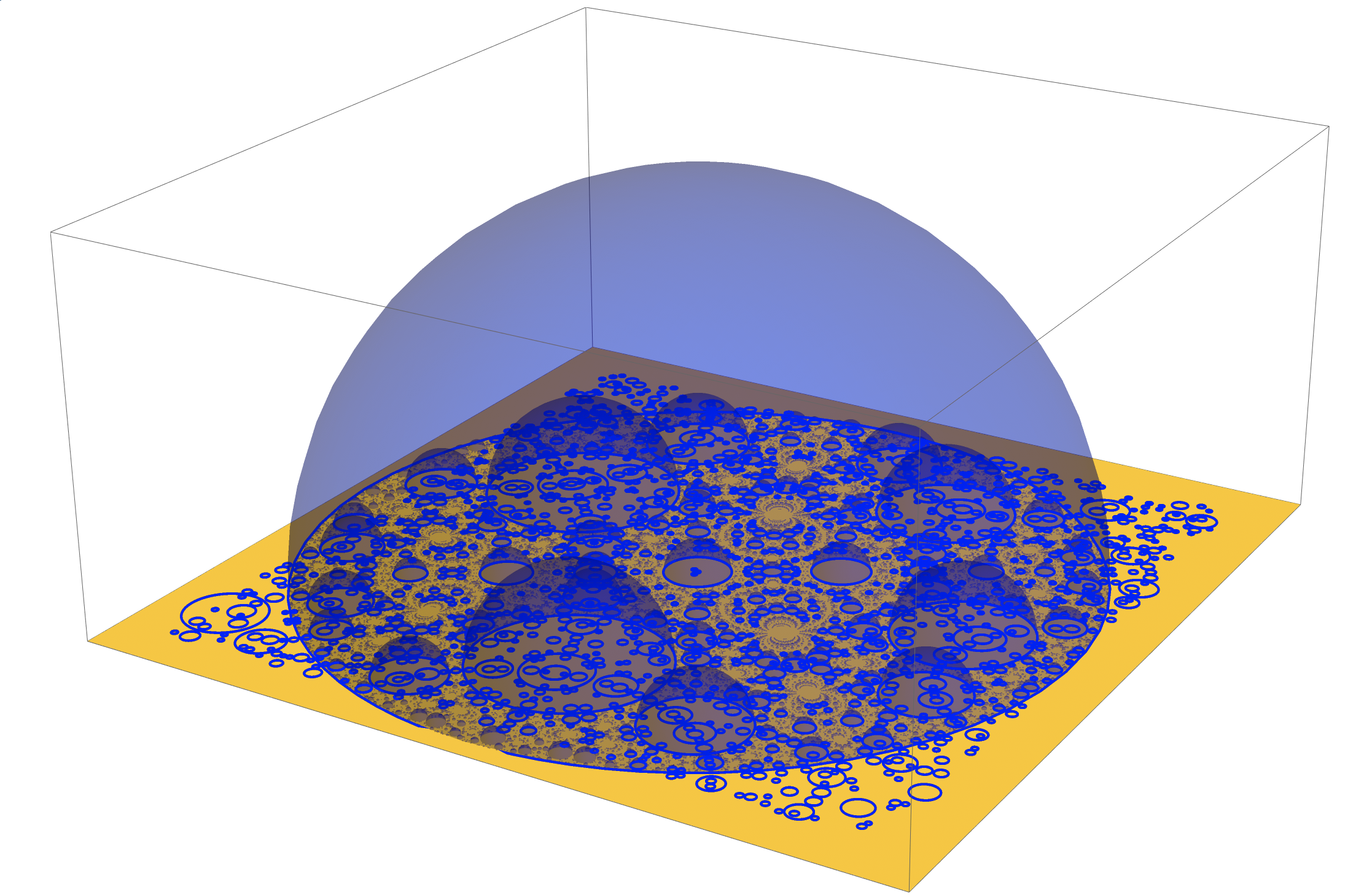}

$(a)$ \hskip1in $(b)$ \hskip1in $(c)$

\caption{$(a)$ A packing $\sP$, $(b)$ its reflection group $\G_\sP$, and $(c)$ its superpacking $\widetilde\sP$.}
\label{fig:superPac}
\end{figure}

Next we %recall %the ``superintegrality'' notion introduced in \cite{KN}. 
%First we 
attach to any sphere packing %(not necessarily crystallographic) 
its ``superpacking.'' 
For an $(n-1)$-sphere $S\subset \bS^n$, denote by $R_S$ reflection through $S$ acting on $\bH^{n+1}$.
Given a packing $\sP$, 
let  
\be\label{eq:GrefDef}
\Ga_{\sP}=\<R_S:S\in\sP\> \ < \ \mathrm{Isom}(\bH^{n+1})
\ee
be the {\bf reflection group of $\sP$}, generated by  
reflections through the spheres in $\sP$. 
The {\bf superpacking} 
%\footnote{$\widetilde\sP$ is independent of the choice of symmetry group $\G$.} 
\be\label{eq:superpack}
\widetilde\sP \ := \ \G_\sP \cdot \sP
\ee
is defined as the orbit of the packing under the action of its reflection group, see \figref{fig:superPac}.
A packing is {\bf superintegral} if its superpacking has all integer bends.\footnote{It turns out that this condition is strictly stronger than integrality; that is, there exist   packings, even crystallographic ones, which are integral but not superintegral \cite{KN}.}
Note that no tangency conditions are imposed on the spheres; indeed the circles in the packing shown in \figref{fig:superPac}$(a)$ are all disjoint (and this packing is superintegral).

As we will show, superintegrality, even in the absence of any other structure imposed on the packing $\sP$, is already %necessarily 
related to (sub)arithme\-ticity, as follows.
Recall that a group of hyperbolic isometries is called ``$k$-arithmetic'' (of simplest type, as assumed throughout) if, possibly after conjugation, it is commensurable with the group 
$O_\cF(\fo)$ 
of $\fo$-integral automorphs of a hyperbolic\footnote{See  \defref{defn:hyperbolic-form}.}  
quadratic form $\cF$
defined over  a totally real number field $k$ with ring of integers $\fo$ (see, e.g., \cite{VS}). 
We call a Zariski-dense, discrete subgroup $\G<\mathrm{Isom}(\bH^{n+1})$ {\bf $k$-subarith\-metic} if $\G$ is contained in a $k$-arithmetic lattice. % group.

\begin{theorem}[Subarithmeticity %SubArithmeticity
Theorem]\label{thm:arithmBugs}
If an orbit 
$$
\cO=\G \cdot S_0
$$
of a fixed sphere $S_0\subset \geo\bH^{n+1}$ under a Zariski dense % discrete 
subgroup $\G<\mathrm{Isom}(\bH^{n+1})$ has all integer bends, then $\G$ is $\Q$-subarithmetic.
More precisely, there exists an isotropic rational hyperbolic quadratic form $\cF$ so that 
$\G$ is contained in $O_\cF(\Z)$.
\end{theorem}

\begin{cor}
In particular, if a general packing $\sP$ %is (somehow) 
happens to be
superintegral, then its reflection group $\G_\sP$ is necessarily  $\Q$-subarithmetic.
\end{cor}

\begin{rmk}
{\em 
These conclusions %do not necessarily 
need not
hold for $\fo$-(super)integral packings, that is, ones with all bends in a ring of integers $\fo$; see \secref{sec:numbFlds}.
}
\end{rmk}

In general, there is not much more one can say about integral or superintegral packings without assuming more structure.

\begin{definition}
A sphere packing $\sP$ is {\bf Kleinian} if its set of limit points also arises as the limit set of a geometrically finite\footnote{
See  \defref{defn:gf}.} 
group $\G_S<\mathrm{Isom}(\bH^{n+1})$.
 We call $\G_S$  a {\bf symmetry group} of the packing.
%As before, a Kleinian packing is {\bf integral} (resp. {\bf superintegral}) if all bends of the packing (resp. superpacking) are integral.
\end{definition}

%We have thus replaced ``finitely generated by reflections'' with ``geometrically finite.'' Recall the definition of the latter: a 

In the special case that $\G$ is generated by finitely many reflections in  hyperplanes,  $\sP$ is called {\bf crystallographic}; such were defined and studied in \cite{KN}.
 It is well-known and  easy  to see (e.g., from the Poincar\'e Fundamental Polyhedron Theorem) that discrete, finitely generated, hyperbolic {\it reflection} groups are geometrically finite, so every crystallographic packing is also Kleinian. 
 
The following statement is both standard and deserves to be stated explicitly.

 \begin{thm}\label{thm:Kleinian-packing} 
Let $\G_S< \mathrm{Isom}(\bH^{n+1})$ be geometrically finite. Then its domain of discontinuity  is a disjoint union of open
%consists of 
round balls (that is, is a Kleinian sphere packing) if and only if the boundary of the convex core $M^*$ of the orbifold $M=\bH^{n+1}/\G_S$ is totally geodesic.
\end{thm}

An important role is played by the ``supergroup'' of a Kleinian packing. If  $\sP$ has symmetry group $\G_S$, then the  {\bf supergroup} $\widetilde\G$ is defined by
$$
\widetilde\G:=\<\G_S,\G_\sP\>.
$$ 
That is, the supergroup is the group generated by both $\G_S$ and $\G_\sP$. Note that both the symmetry group $\G_S$ and supergroup $\widetilde\G$ are not uniquely determined by the packing $\sP$; indeed, any nontrivial normal subgroup of $\G_S$ will have the same limit set. 
{\it A priori} it is not even obvious that $\widetilde\G$ acts discretely, but  it is in fact a lattice, acting on $\bH^{n+1}$ with finite covolume; see the Structure \thmref{thm:lattice}. 

\medskip 
We turn our attention now to the integral and superintegral Kleinian packings. 
An immediate corollary of the Subarithmeticity \thmref{thm:arithmBugs} is that, if a Kleinian packing $\sP$ is integral, then 
 any symmetry group $\G_S$ is $\Q$-subarithmetic.
If $\sP$ is moreover superintegral, 
then 
its supergroup $\widetilde\G$ is itself $\Q$-arithmetic.
This is because the superpacking $\widetilde\sP$ can also be given as the orbit of the packing $\sP$ under the action of $\widetilde\G$, together with the obvious fact that,
if a group is subarithmetic and a lattice, then it is arithmetic!

In the case of crystallographic packings, the main result of \cite[Thm 18]{KN} was the following 
finiteness theorem. Before stating the theorem, note that it is shown in \cite[Thm 3]{KN} that  there exist infinitely many {\it conformally inequivalent} 
superintegral crystallographic packings in certain dimensions up to $n=18$.  We say that two Kleinian packings are {\bf commensurable} if (conjugates of) their supergroups are.

\begin{thm}[Finiteness Theorem \cite{KN}]\label{thm:fin}
Superintegral crystallogra\-phic packings exist in only finitely-many dimensions, and there are finitely many in each dimension, up to commensurability. 
\end{thm}

Indeed, if a crystallographic packing is superintegral, then its supergroup is a  $\Q$-arithmetic {\it reflective} lattice.
Arithmetic reflective lattices are known (see, e.g., the discussion in \cite{Bel}) to lie in finitely many commensurability classes in finitely many dimensions, which implies  \thmref{thm:fin}.
Moreover, there are no $\Q$-arithmetic reflective lattices acting on $\bH^{n+1}$ with $n+1=20$ or $n+1\ge22$, see \cite{Ess}, so there are no superintegral crystallographic packings in $n=19$ or $n\ge21$ dimensions. Worse yet, crystallographic packings are not yet known to exist (nevermind integrality) in dimensions $n=14, 15, 16, 18$ and $20$, although reflective lattices are known in one more than these dimensions \cite{KN}.

\medskip
While the classification of commensurability types of superintegral crystallographic packings awaits first that of arithmetic hyperbolic reflection groups, one can  completely  classify superintegral Kleinian packings in terms of arithmetic groups. 
The following is the first main theorem of this paper.

\begin{theorem}[Classification Theorem]\label{thm:class}
A hyperbolic lattice is commensurable to a %the 
supergroup of
 a superintegral Kleinian packing
 if and only if
it is a non-uniform $\Q$-arithmetic lattice of  simplest type.
 \end{theorem}

So in contradistinction with the Finiteness \thmref{thm:fin} for crystallographic packings, more general Kleinian packings exist in every dimension. Moreover, the Classification \thmref{thm:class} answers a question posed in \cite{KN} on whether there exist superintegral crystallographic packings with {\it cocompact} supergroups: there do not!
We emphasize again that the spheres in a superintegral Kleinian (or even just crystallographic) packing could all be disjoint, but the supergroup $\widetilde\G$ must have cusps. See \exref{ex:nonTangent}.

\begin{rmk}\label{rmk:10}
{\em
The precise role of $\Q$-isotropy of the corresponding quadratic form leading to (super)integrality of packings is elucidated in \lemref{lem:bbIs}, which shows that the ``covector'' corresponding to the bend is itself isotropic in the dual form.
}
\end{rmk}

\begin{remark}\label{rmk:nonInt}{\em
Note that even if
a symmetry group $\G_S$ is geometrically finite and
 the convex 
core $M^*$ of $M=\bH^{n+1}/\G_S$ has totally geodesic boundary (so its limit set  gives rise to 
 a Kleinian packing $\sP$), 
and if the supergroup $\widetilde\G$ of $\sP$ is non-uniform and $\Q$-arithmetic of simplest type,
it still need {\it not} be the case that the packing $\sP$ is necessarily integral, for any conformal
choice of coordinates on $\bS^n$, see  \exref{ex:nonInt}. 
But the next theorem states that, if  $M^*$  has only one boundary component (that is, $\G_S$ acts transitively on the spheres in the packing), 
then the packing {\it is} necessarily superintegral.
%see also \rmkref{rmk:intNonInt}.
}
\end{remark}

\begin{thm}\label{thm:converseArith}
$(i)$ Suppose that a supergroup $\widetilde\G$ of a packing $\sP$ is a non-uniform $\Q$-arithmetic of simplest type (i.e.  is 
commensurable to $O_\cF(\Z)$, where $\cF$ is a rational hyperbolic quadratic form) and, moreover, a symmetry group $\G_S< \widetilde\G$ acts transitively on the  set of
spheres in the packing $\sP$. Then there is always a conformal change of coordinates such that $\sP$ is superintegral. 

$(ii)$ More generally, let $\G$ be any non-uniform $\Q$-arithmetic group of the simplest type and 
 let $S$ be a sphere for which the reflection $R_S$ lies in $\G$.  Then the $\G$-orbit
 of $S$ can be made integral by a suitable conformal change of coordinates.
\end{thm}

There are two main ingredients leading  to the proof of the Classification \thmref{thm:class}, one geometric (constructing a Kleinian packing from the Structure \thmref{thm:lattice} below), and the second arithmetic, to ensure that the packing thus constructed is indeed superintegral; it is the second stage that forces $\Q$-arithmeticity and non-uniformity of the lattice.

A key  step in the geometric argument
relies on  Millson's theorem \cite{Millson} that every arithmetic hyperbolic lattice of simplest type is commensurable to a lattice $\widetilde\Ga$ such that $\H^{n+1}/\widetilde\Ga$ contains a nonseparating totally geodesic complete 
%co-dimension-1
hypersurface of finite volume.

We will also show that, in each commensurability class, there are not only infinitely many conformally inequivalent packings (as shown for crystallographic packings in \cite{KN}), but their Hausdorff dimensions can be made arbitrarily close to maximal.

\begin{theorem}[Abundance Theorem]\label{thm:abundance}
In every dimension $n\ge2$, there exist quasiconformally inequivalent  superintegral Kleinian packings, whose limit sets have Hausdorff dimension approaching $n$. 
%When $n=2, 3$, these packings can be taken by pairwise disjoint balls and, thus,  topologically-equivalent.}
\end{theorem}

%The issue with topological equivalence is unclear: We get it for sure, \cite{DV} or \cite{Cannon}, if the spheres in packing are %pairwise disjoint. But our super-lattice 
%has to be non-uniform. Thus, we need to find non-uniform lattices which contain uniform codimension 1 sublattices. 
%This can be done in low dimensions....

Note the contrast to Phillips-Sarnak's and Doyle's work on Schottky groups \cite[Theorem 5.4]{PhillipsSarnak1985}, \cite{Doyle1988}, showing that their limit sets have Hausdorff dimensions bounded strictly away from the ambient dimension $n$.
Also note that the conformally inequivalent superintegral Polyhedral Circle Packings in $\geo\bH^3$ constructed in \cite[Thm 7]{KN} are all Schottky, so their dimensions cannot approach $n=2$.

\subsection{Packings over Number Fields}\label{sec:numbFlds}\

To complete the discussion of Kleinian packings (before we turn to Kleinian ``bugs''), we show that the theory breaks down over number fields, answering another question posed in \cite{KN}. For a totally real number field $k$ with ring of integers $\fo$, a packing is $\fo$-integral (resp. $\fo$-superintegral) if every bend in its packing (resp. superpacking) lies in $\fo$.

\begin{prop}\label{prop:golden}
The Subarithmeticity Theorem does not hold in number fields.
\end{prop}

Indeed, already for the simplest case of the golden ring $k=\Q(\sqrt 5)$, there exist $\fo$-superintegral crystallographic packings whose supergroups are non-arithmetic. (See \exref{ex:golden}.)
But this is not the whole story, as we also have the following.

\begin{theorem}\label{thm:quadExt}
For any $k$-arithmetic hyperbolic lattice of simplest type $O_\cF(\fo)$,
%  $\widetilde\G<\mathrm{Isom}(\bH^{n+1})$, which is a $k$-arithmetic lattice of simplest type,
there is a totally real, quadratic extension $k'\supset k$, with ring of integers $\fo'$, such that
$O_\cF(\fo)$ %$\widetilde\G$
 is commensurable to the supergroup of 
 an $\fo'$-superintegral Kleinian packing. 
 \end{theorem}

Note that there is no condition on being non-uniform here, unlike \thmref{thm:class}. Thus, even the $\Q$-anisotropic form $x_1^2+x_2^2+x_3^2-7x_4^2$ has  $\fo'$-superintegral Kleinian packings, for certain quadratic rings $\fo'$.
This leaves  open the problem of giving a proper formulation for how to characterize which $\fo'$-superintegral packings come from arithmetic groups.

\subsection{Kleinian Bugs}\

To further extend the notion of a packing, we will allow spheres to meet at a finite set of dihedral angles, as follows.

\begin{figure}
\begin{center}

\begin{tikzpicture}[node distance=1cm, auto,]
 %nodes
 \node (market) {\includegraphics[width=.7\textwidth]{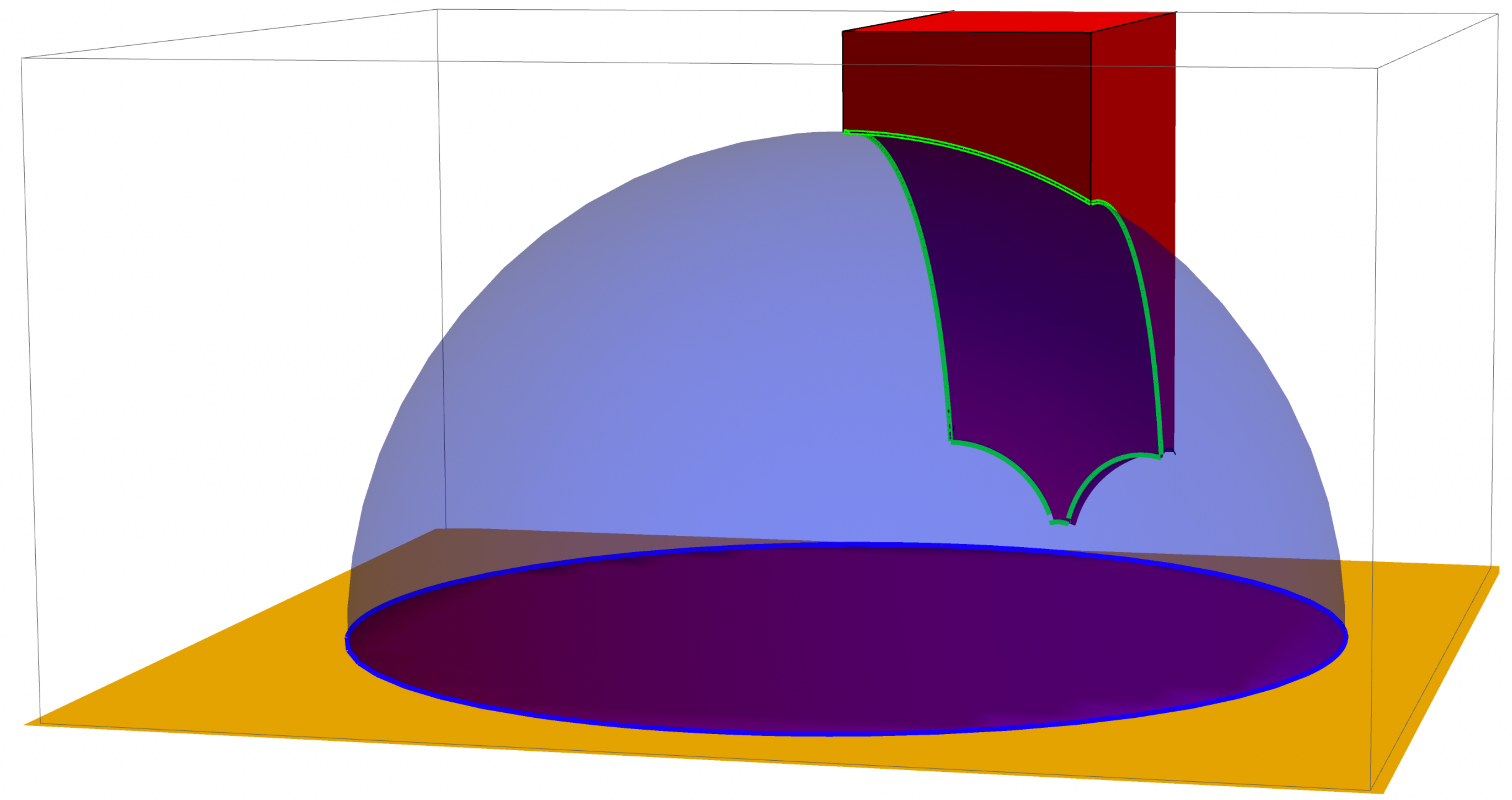}};
 \node at (1,2) (Pend) {};
 \node at (-2,3) (Pstart) {$P$}
     edge[pil,bend left=30] (Pend);
 \node at (1.25,0.5) (Fend) {};
 \node at (-4,2.5) (Fstart) {$F_\alpha$}
     edge[pil,bend right=45] (Fend);
 \node at (-2,0) (Halend) {};
 \node at (-5,1.5) (Halstart) {$H_\alpha$}
     edge[pil,bend right=10] (Halend);     
\node at (5.25,-.3) (ddHn) {$\partial_\infty \H^{n+1}$};
 \node at (-1,-1.3) (BHalend) {};
 \node at (-5,-1) (BHalstart) {$B_\alpha$}
     edge[pil,bend left=30] (BHalend);     
 \node at (3.4,-1.5) (ddHalend) {};
 \node at (5,-2) (ddHalstart) {$\partial_\infty{H_\alpha}$}
     edge[pil,bend right=45] (ddHalend);     
\end{tikzpicture}

\end{center}
\caption{A portion of a polyhedron $P$ having bounding facet $F_\alpha$ lying in the hyperplane $H_\alpha$. Also shown is the 
 ideal boundary 
 $B_\alpha\subset \partial_\infty\H^{n+1}$
 of the  complementary half-space $H_\alpha^-$,
 and the 
 cooriented round sphere $\partial_\infty H_\alpha$
 bounding  $B_\alpha$.
}
\label{fig:Hga}
\end{figure}
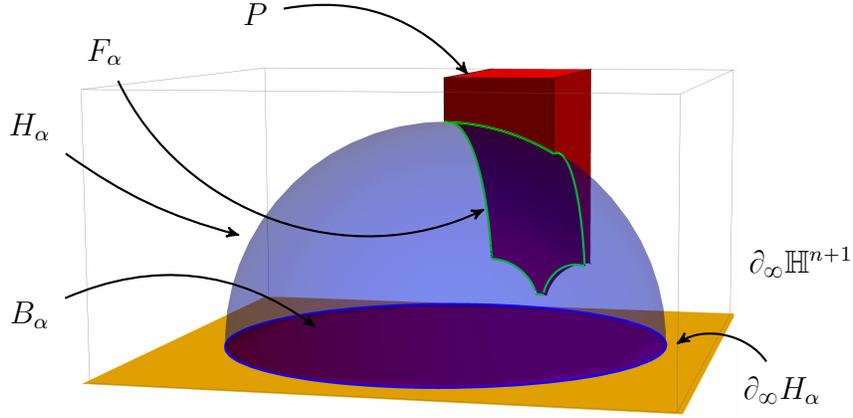

A {\em cooriented} round sphere in $\bS^n$ is a round sphere together with a choice of a nowhere vanishing normal vector field. 
Cooriented spheres are in a natural bijective correspondence with cooriented hyperbolic hyperplanes in $\H^{n+1}$. A {\em convex polyhedron}, $P$, in $\H^{n+1}$ is the intersection of a locally finite family of hyperbolic half-spaces $H_\alpha^+$ 
in $\H^{n+1}$. It suffices to consider only half-spaces $H_\alpha^+$ such that the hyperplanes $H_\alpha$ intersect $\partial P$ along  facets $F_\al$. 
We coorient these hyperplanes $H_\alpha$  so that the normal vector fields point into the complementary half-spaces  
$H_\alpha^-$. The ideal boundary $\geo H_\al^-$ of 
the half-space $H_\al^-$ is a round ball $B_\al$ 
bounded by the cooriented round sphere $\geo H_\al$, see \figref{fig:Hga}.

\begin{definition}
A {\bf bug} $\sB$ is an infinite collection of cooriented spheres in $\bS^n\cong\geo\bH^{n+1}$ containing facets of a convex polyhedron $P=P_\sB\subset \H^{n+1}$, with dihedral  
 angles lying in a finite subset of 
$$
 \frac\pi\N \ = \ \left\{\frac\pi m\ : \ m\in\N\right\},
$$
such that the union of round balls bounded by these spheres is dense in $\bS^n$.  
As before, a bug is {\bf integral} if its spheres all have integer bends.
\end{definition}

In plainer terms (forgetting coorientations), a bug is a collection of spheres, any pair of which is either disjoint, tangent or which intersect at angles $\pi/m$, for finitely many $m\in\N$.
See \figref{fig:1} for an integral bug, and  \figref{fig:2} for the corresponding hyperbolic polyhedron $P$, which should justify our calling these objects ``bugs.'' 

\subsubsection{Geometric Aspects of  Bugs}\

Before turning to any arithmetic properties of bugs, we discuss purely geometric aspects, culminating in the Structure \thmref{thm:lattice}.
As in \eqref{eq:GrefDef}, the {\em reflection group}
%\footnote{The reader confused by the fact that the reflection group is denoted by $\Ga_\sB$ rather than $\Ga_R$ might be relieved by learning that the Russian letter ${\mathrm P}$ corresponds to the English letter $R$! 
%{\color{red}Ha! Nice obvservation! 
%Though we should probably leave it out for the paper...
%}} 
 generated by reflections through the spheres in $\sB$, denoted $\Ga_{\sB}%=\Ga_\sB
 $,  is discrete and  has as its fundamental domain 
the %a %%% Not "a", it's   
 hyperbolic convex polyhedron 
$P=P_\sB$ bounded by the hemispheres.

\begin{figure}
\includegraphics[width=.45\textwidth]{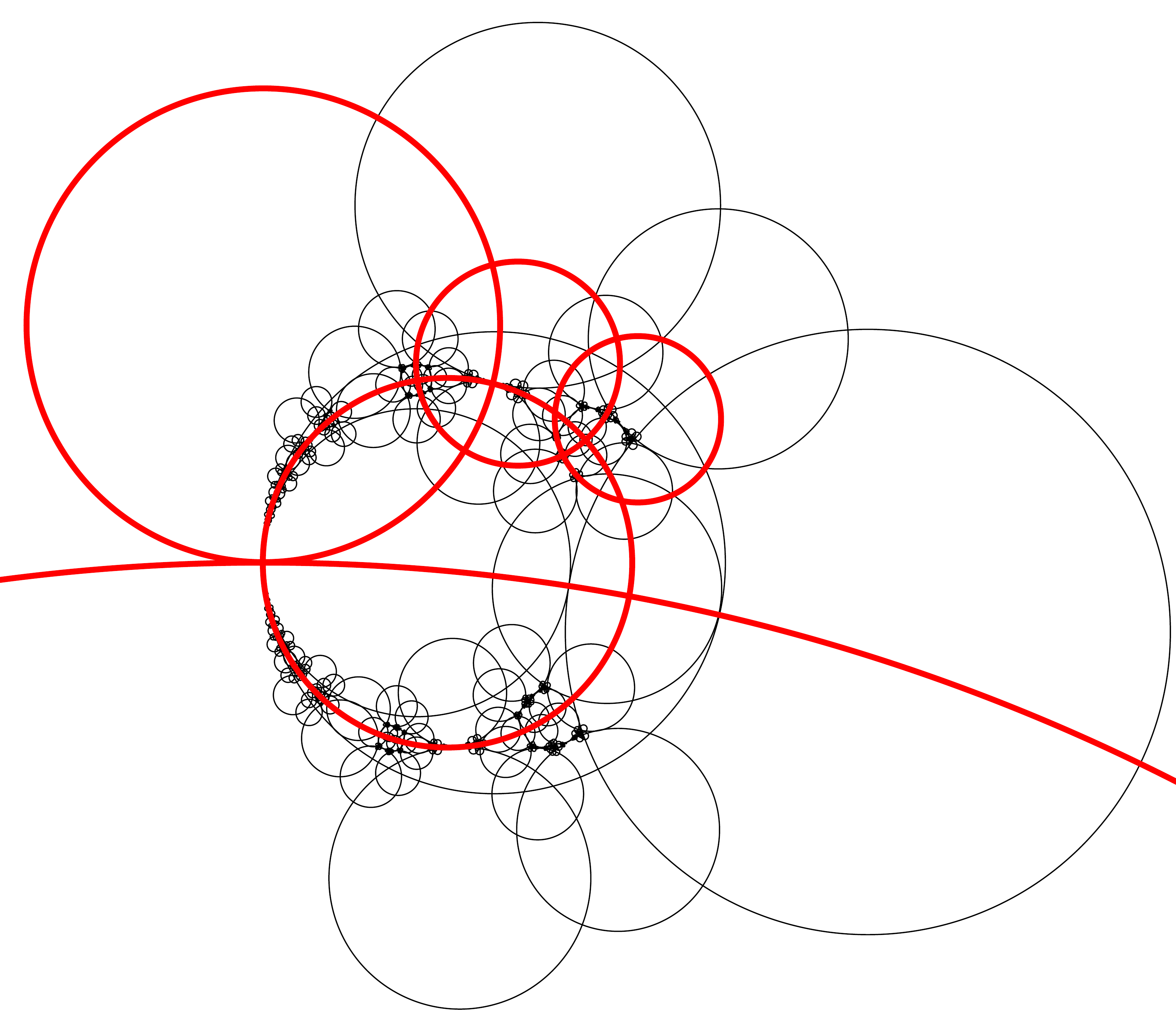}\hskip.3in
\includegraphics[width=.45\textwidth]{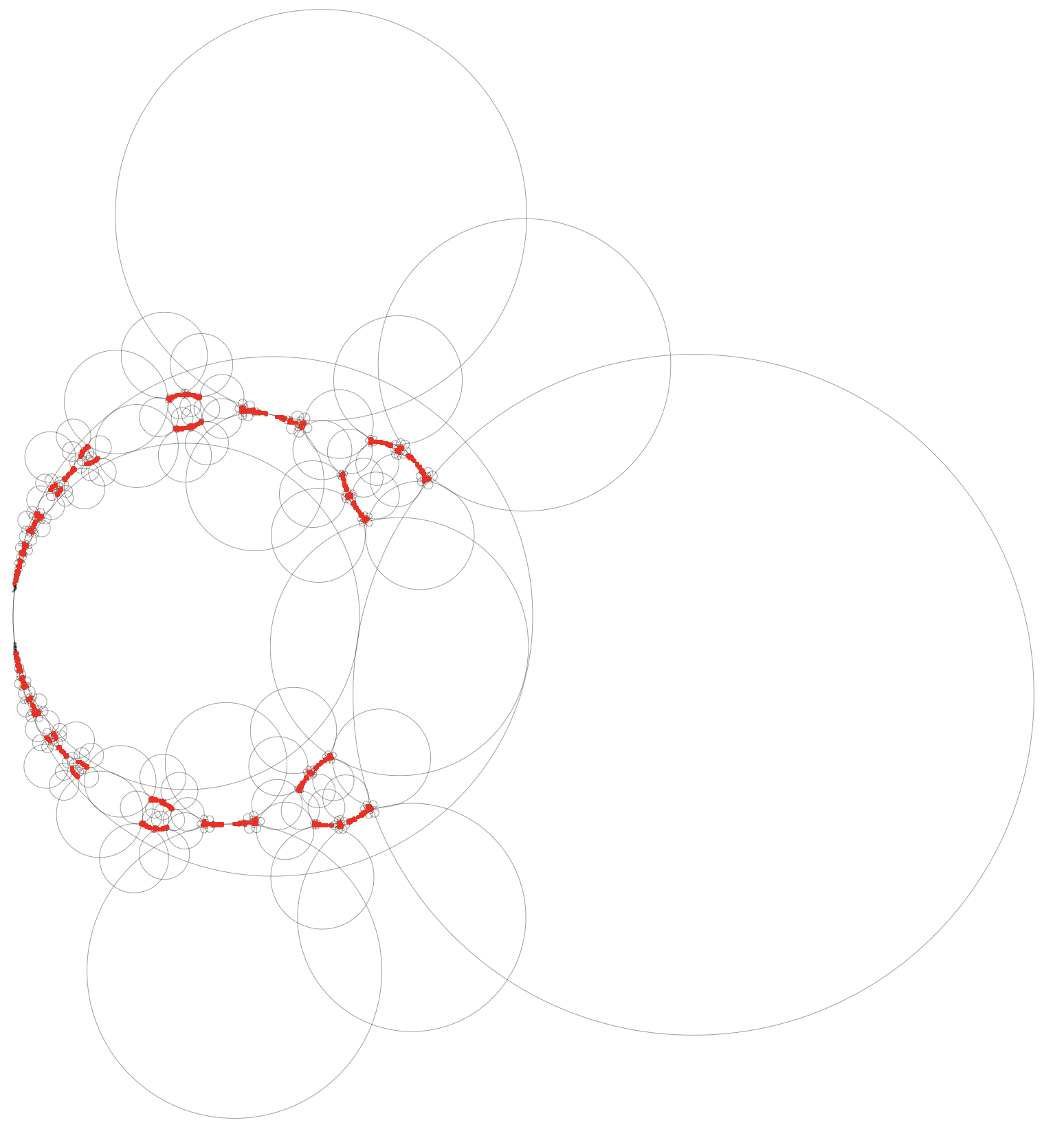}

$(a)$\hskip1in$(b)$
\caption{$(a)$ The bug is Kleinian with symmetry group generated by reflections in the thick red circles. $(b)$ The accumulation set in $\geo\bH^{n+1}$ does not in general contain spheres.}
\label{fig:3}
\end{figure}

In this setting, we need a substitute for the limit set of the packing, which will no longer necessarily contain spheres, see \figref{fig:3}(b).

\begin{definition}
The {\bf accumulation set} $\A(\sB)$ of a bug $\sB$ 
consists of those points $\xi\in \dd_\infty\H^{n+1}$ such that every neighborhood of $\xi$  contains infinitely many spheres in $\sB$.
\end{definition}

As a first (it turns out, insufficient) step, we define ``geometric'' bugs in terms of symmetry groups as before:

\begin{figure}
\includegraphics[width=.75\textwidth]{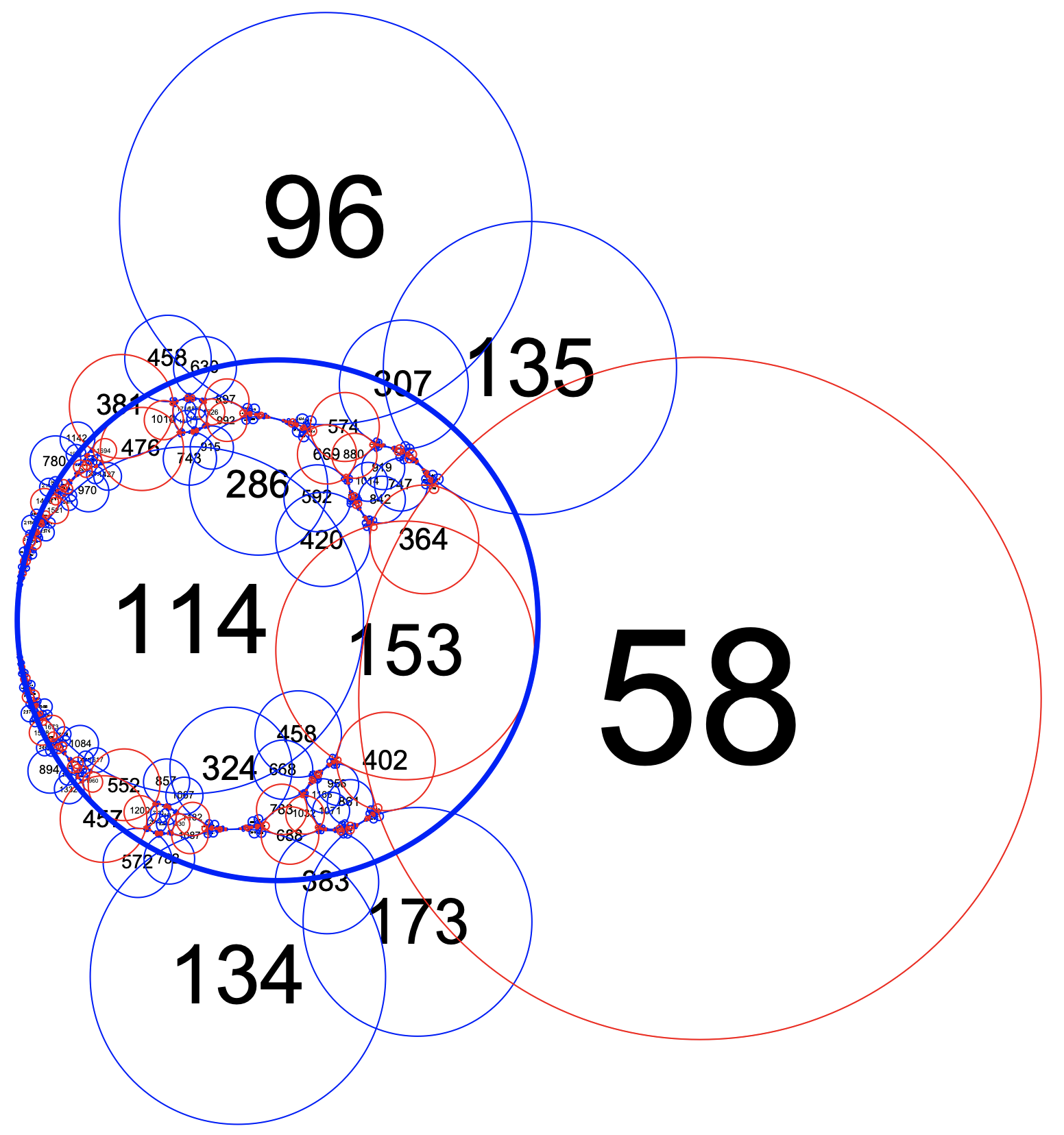}\
\caption{An integral bug with all dihedral angles $\pi/3$. A number at the center of a circle denotes its bend, that is, inverse radius. The thick unlabeled blue circle has bend $-76$. The color scheme is as described in \rmkref{rmk:constr}.}
\label{fig:1}
\end{figure}

\begin{figure}
\includegraphics[width=.65\textwidth]{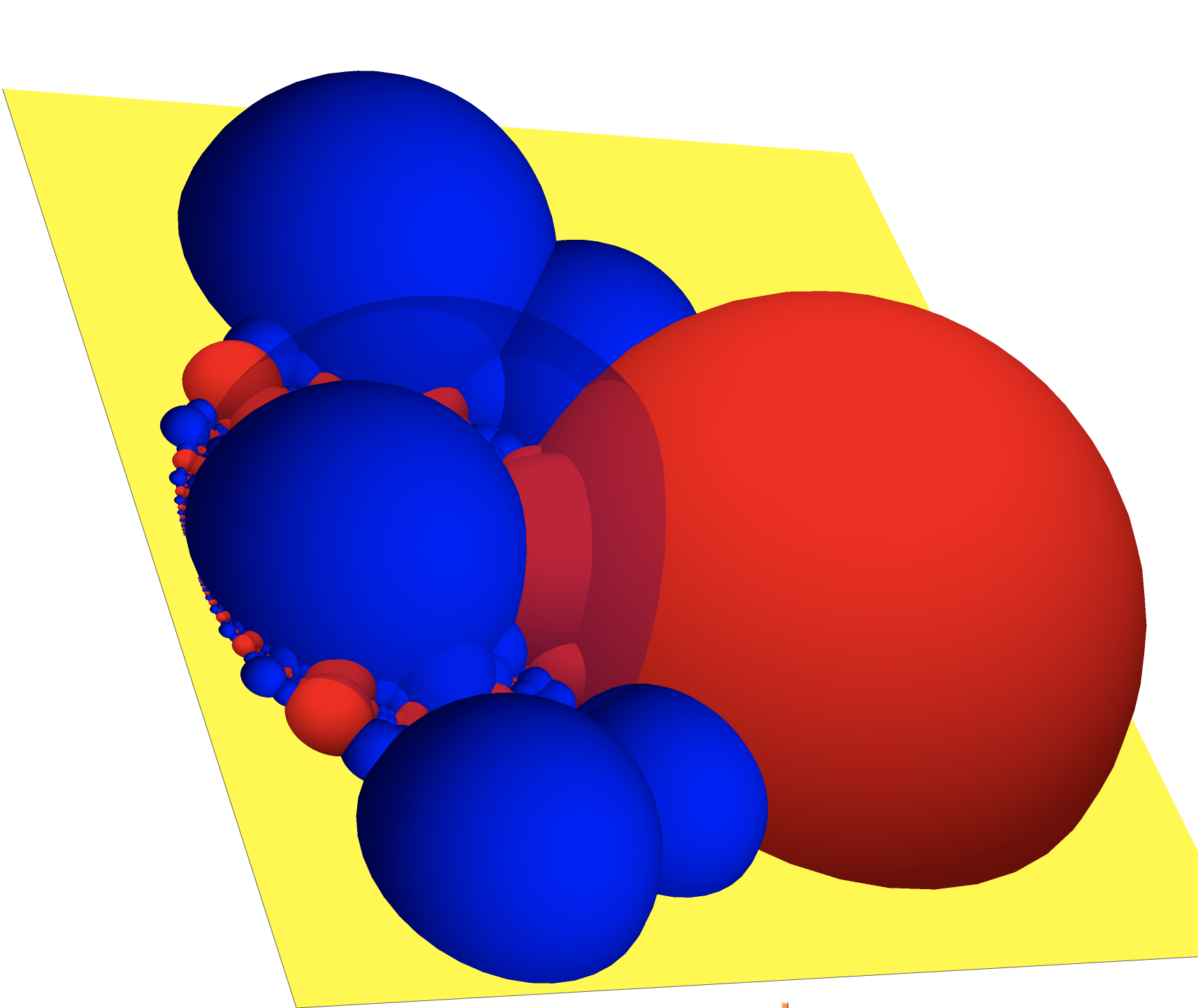}
\caption{View from $\bH^3$ of the bug polyhedron $P$  from \figref{fig:1}. The polyhedron is exterior to all solid hemispheres shown, and interior to the one opaque hemisphere, which has bend $-76$.}
\label{fig:2}
\end{figure}

\begin{definition}
A bug $\sB$ is {\bf geometric} if there exists  a geometrically finite subgroup $\G_S<\mathrm{Isom}(\bH^{n+1})$ such that
\begin{enumerate}
\item 
the accumulation set of the bug
is the limit set of $\G_S$, and
\item 
 $\G_S$ preserves the bug polyhedron $P=P_\sB$. 
\end{enumerate}
Such a $\G_S$ is called a {\bf symmetry group} of the bug.
\end{definition}

Recall  again that a symmetry group of a bug is not unique. But its commensurability class is uniquely determined by the bug, provided that the accumulation set of ${\sB}$ 
consists of more than two points; see \lemref{lem:commensurable}. 

The supergroup is defined in the now-familiar way:

\begin{definition}
Given a bug $\sB$ and a symmetry group $\G_S$, the {\bf supergroup} 
$$
\widetilde\G \ := \ \<\G_S, \Ga_\sB\> \ < \ \mathrm{Isom}(\bH^{n+1}),
$$
is defined
as the group generated by the symmetry group $\G_S$ and the reflection group $\Ga_\sB$.
\end{definition}

Naively, one could have expected that, as in the crystallographic or Kleinian packing setting, the supergroup of a bug is a lattice; unfortunately this is false in general (see \exref{ex:nonCusp}), owing to the possible ``incompatibility'' of parabolic subgroups.
Recall that a discrete subgroup $\Pi<\mathrm{Isom}(\H^{n+1})$ is called {\em parabolic}\footnote{Sadly, this classical terminology is inconsistent with the one used in the theory of algebraic groups.}  if it contains a parabolic element $g$ and stabilizes the fixed point of $g$.  Such a subgroup consists of parabolic and elliptic elements. 

\begin{definition}\label{def:cuspComp}
An ordered pair of  subgroups $\Ga_1, \Ga_2< \mathrm{Isom}(\bH^{n+1})$ is {\bf cusp-compatible} if, for every maximal parabolic subgroup $\Pi_2< \Gamma_2$, there exists a subgroup  $\Pi_1< \Gamma_1$  such that the subgroup $\Pi$ generated by $\Pi_1, \Pi_2$ is virtually abelian of virtual rank $n$. 
\end{definition}

Note that cusp compatibility depends only on the commensurability classes of $\Ga_1, \Ga_2$. 

\begin{definition}
A bug $\sB$ is called {\bf Kleinian} if it is geometric with symmetry group $\G_S$, such that the ordered pair $\Ga_\sB$, $\G_S$ is cusp-compatible.
\end{definition}

We remark that, if $\sB$ is not just a bug but also a packing, then the cusp-compatibility condition is automatic, see \lemref{lem:simple->lattice}.

In \secref{sec:lattice} we will prove the second main theorem of this paper, which is the following extension of the Crystallographic Structure Theorem in \cite[Thm 28]{KN}.

\begin{thm}[Structure Theorem for Kleinian Bugs and Packings]\label{thm:lattice}
1. If $\sB$ is a Kleinian bug, then its supergroup $\widetilde\G$ is a lattice.

2. (i) Conversely, suppose we are given a lattice $\widetilde\G<\mathrm{Isom}(\bH^{n+1})$,  a convex fundamental polyhedron 
%domain 
$D$ of  $\widetilde\G$, and a finite set $S'$ of elements $\g_j\in\widetilde\G$ pairing the facets %bounding walls 
of $D$.
%; that is, $S'$ is a minimal set of generators for $\widetilde\G$. 
%Although not all minimal generating sets appear this way 
Let $R\subset S'$ be a nonempty subset consisting of reflections (assuming such exist) in some facets $F_\al\subset \partial D$, with $\al$ in an indexing set $A$.  As before, let $H_\al$ be the cooriented hyperbolic 
hyperplanes containing the facets $F_\al$. 
Let $S:=S'\setminus R$ and $\G_S:= \<S\><\widetilde\G$. Then 
the orbit of the set of cooriented spheres in $ \{\geo H_\al; \al\in A\}$ under the group $\G_S$,  
$$
\sB \ := \ \G_S \cdot \{\geo H_\al; \al\in A\},
$$
is a Kleinian bug with a symmetry group $\Ga_S$. %{\color{red} Maybe there's a cleaner way to say this?}

(ii) Suppose, moreover,  that: (a) the facets $F_\al$, $\al\in A$, have pairwise disjoint or tangent cooriented spheres $\geo H_\al$, (b) the hyperplanes $H_\al$ meet the other bounding walls of $D$ either tangentially or orthogonally (or not at all).  
Then $\sB$ is in fact a Kleinian {\bf packing}, that is, if two spheres of $\sB$ intersect, they do so tangentially.
\end{thm}

In particular, every lattice $\widetilde\G<\mathrm{Isom}(\bH^{n+1})$ containing a reflection yields a Kleinian bug. 

\begin{rmk}{\em
While it is not hard to show that the supergroup of a crystallographic or Kleinian packing 
is a lattice, the proof of \thmref{thm:lattice} is rather more involved, see \secref{sec:lattice}.
}
\end{rmk}

\begin{rmk}\label{rmk:constr}
{\em
We applied the above construction to create the bug in Figures \ref{fig:1}--\ref{fig:2} from the extended Bianchi group $\widehat{Bi}(19)$\footnote{\label{foot:Bianchi}Here and throughout, $\widehat{Bi}(D)$ denotes the extended Bianchi group, that is, the maximal discrete subgroup of $\Isom(\H^3)$  containing the group $\PSL_2(\cO_{D})$, where $\cO_{D}$ is the ring of integers of $\Q(\sqrt{-D})$.};
% is not reflective (see \cite{BelolipetskyMcleod2013}); 
 namely, two reflective walls were chosen for $R
\subset S'$ from this lattice,
%\footnote{The group we chose was the reflection subgroup of the maximal discrete extension of the Bianchi group $\SL_2(\cO_{19})$, with $\cO_{19}$ the ring of integers of $\Q(\sqrt{-19})$.}, 
and the bug was created as the orbit of these two walls under the remaining generators. The orbit of one wall was colored blue and the other red.
%The code is available at \url{asfd}.
}
\end{rmk}

\subsubsection{Arithmetic Aspects of  Bugs}\

Turning our attention to arithmetic properties of bugs, we mimic \eqref{eq:superpack} with the following.

\begin{definition}
The {\bf superbug} $\widetilde\sB$ of an (arbitrary) bug $\sB$
is defined as the orbit of the bug under its reflection group, $\G_\sB$:
$$
\widetilde\sB\ : = \ \G_\sB \cdot \sB.
$$
A bug is {\bf superintegral} if every sphere in its superbug has integer bend.
\end{definition}

The Subarithmeticity \thmref{thm:arithmBugs} applies just as well to bugs, giving the following:

\begin{thm} \
\begin{itemize}
\item
If a bug (with no further structure imposed) is superintegral, then its reflection group $\G_\sB$ is $\Q$-subarithmetic.
\item
If a bug is Kleinian and integral, then its symmetry group $\G_S$ is $\Q$-subarithmetic.
\item
If a bug is Kleinian and superintegral, then its supergroup $\widetilde\G$ is $\Q$-arithmetic.
\end{itemize}
\end{thm}

In the opposite direction, we already have from the Classification \thmref{thm:class} that every non-uniform $\Q$-arithmetic lattice (of simplest type) is commensurable to the supergroup of a superintegral Kleinian packing (and hence bug). Similarly, we have from  \thmref{thm:converseArith} that if a symmetry group $\G_S$ (resp. supergroup $\widetilde\G$) of a Kleinian bug  is non-uniform and $\Q$-subarithmetic (resp. $\Q$-arithmetic), then there is a conformal choice of coordinates for which the bug can be made integral (resp. superintegral).

\begin{rmk}{\em
When a Kleinian bug (or just packing) is integral but not superintegral, we have the following curious situation: its symmetry group $\G_S$  is a subgroup of some integer orthogonal group $O_{\cF}(\Z)$,} but
 it  is also contained in a (non-arithmetic) lattice, the supergroup $\widetilde\G$.
\end{rmk}

Note further that there are essentially only three ways that a lattice $\widetilde\G<\mathrm{Isom}(\bH^{n+1})$ (with $n\ge2$) can be non-arithmetic. By %Mostow/
Selberg's Rigidity \cite{Sel}, after conjugation, $\widetilde\G$ is contained in some $O_{n+1,1}(k)$ where
 %$q$ is a quadratic form of signature $(n+1,1)$ over a 
 $k$ is a
 number field %$k$ (
 with a given embedding to $\R$. For $\widetilde\G$ to be non-arithmetic, either: 
 \begin{itemize}
 \item[$(i)$] $k$ is not totally real, or
 \item[$(ii)$] $k$ is totally real but for at least one non-identity embedding $\iota$, the orthogonal group $O_{n+1,1}(\iota (k))$ is non-compact (that is, the quadratic from is not hyperbolic), %quadratic form $q$ is indefinite, 
 or 
 \item[$(iii)$] $k$ is totally real and the quadratic form is hyperbolic, 
 %(that is, has only the one indefinite embedding), 
 but  the entries of $\widetilde\G$ as elements of $k$ have  unbounded denominators. 
 \end{itemize}
Lattices satisfying the last condition $(iii)$ are called {\it quasi-arithmetic} by Vinberg.
The only currently known integral but non-superintegral %packings 
bugs (e.g. the crystallographic packings constructed in \cite[Lemma 20]{KN})
all come from this last situation (with $k=\Q$).
It is interesting to investigate whether this is the only non-arithmeticity type possible for integral but not superintegral bugs.
(Note that for $\fo$-superintegrality, the non-arithmeticity in \propref{prop:golden} is of type $(ii)$; see \exref{ex:golden}.)

\subsection{Outline of the paper}\

In \secref{sec:lattice} we prove geometric results on bugs, including the key the Structure \thmref{thm:lattice}. We spend \secref{sec:arith} discussing the arithmetic properties of bugs, including the main Subarithmeticity  \thmref{thm:arithmBugs}.
Of particular interest to the reader may be \secref{sec:exs} where a number of archetypal examples are constructed.
 We conclude in \S\ref{sec:Hdim} with a discussion of Hausdorff dimensions of accumulation sets of superintegral bugs, proving the Abundance \thmref{thm:abundance}.
%Finally, we use the \appref{sec:comb} to discuss methods for making the above constructions {\it finitary}, that is, not checking equality of infinite limit and accumulation sets.

%\subsection{Notation}
\subsection*{Acknowledgements}\

The second-named author  would like to thank Curt McMullen and Peter Sarnak for many enlightening conversations and suggestions.

%\newpage

\section{Proof of the Structure \thmref{thm:lattice}}\label{sec:lattice}

In what follows, $d$ denotes the hyperbolic distance on $ \H^{n+1}$. For a subset $A\subset \H^{n+1}$ the {\em ideal boundary} of $A$, denoted $\geo A$,  
is the accumulation set of $A$ in the ideal boundary sphere $\bS^n=\geo \H^{n+1}$ of $\H^{n+1}$. 

Recall that for a discrete subgroup $\Ga< \Isom(\bH^{n+1})$, the {\em convex hull} of the limit set $\Lambda$ of $\Ga$, denoted $C=C_\Ga$, is the intersection of all $\Ga$-invariant closed convex nonempty subsets in $\bH^{n+1}$. The convex hull has the property that $\geo C=\Lambda$ unless $\Lambda$ is a singleton (in which case $C=\emptyset$). The 
{\em convex core} $M^*$ of the orbifold $M=\bH^{n+1}/\Ga$ is the quotient $C/\Ga$.  

\begin{definition}\label{defn:gf}
A discrete  subgroup $\G<\mathrm{Isom}(\bH^{n+1})$ is geometrically finite if $\G$ is virtually torsion free and there is an $\vep$-thickening of 
 the convex core $M^*$ of the orbifold $M=\bH^{n+1}/\Ga$
that has finite volume.
\end{definition}

A sufficient condition for geometric finiteness is the existence of a finitely-sided fundamental polyhedron of $\G$ in $\H^{n+1}$. 

\begin{definition}
A discrete subgroup $\Ga< \Isom(\bH^{n+1})$ is called {\em convex-cocompact} if its limit set is not a singleton and 
the convex core $M^*$ 
of the orbifold $M=\bH^{n+1}/\Ga$ is compact. Equivalently, there exists a closed convex nonempty $\Ga$-invariant subset of 
$\bH^{n+1}$ on which $\Ga$ acts cocompactly. 
\end{definition}

Convex-cocompact subgroups have the property that 
\begin{equation}\label{eq:inrad} 
\inf \{ d(gx, x): g\in \Ga^*, x\in C_\Ga\}>0,
\end{equation}
where $\Ga^*\subset \Ga$ consists of elements of infinite order. 

We refer to \cite{Bowditch93, Rat} for more background on geometrically finite and convex-cocompact isometry groups of hyperbolic space.

\medskip 
{\bf Part 1.} Let $\sB$ be a Kleinian bug and let $P\subset \H^{n+1}$ be the corresponding convex polyhedron. We wish to show that the supergroup $\widetilde\G$ is a lattice.  

If $\Phi\subset \H^{n+1}$ is a convex fundamental domain of an arbitrary  discrete subgroup $\Ga < \Isom(\H^{n+1})$, 
then the $\Ga$-orbit of the relative interior of $\geo \Phi$ (with respect to $\bS^n$) is dense in the domain of discontinuity of $\Ga$ in $\bS^n$. 
Since $\sB$ is a bug, for $\Phi=P$, this relative interior of $\geo P$ is empty and, hence the limit set of $\Ga_{\sB}$ is equal to the entire sphere $\bS^{n}$.

The ideal boundary $\geo P\subset \geo \H^{n+1} = \bS^{n}$ is the disjoint union 
\be\label{eq:geoP}
\geo P = \A(P)\sqcup \V(P),
\ee 
where each point of  $\V(P)$ is isolated in the accumulation set $\A(P):=\A(\sB)$ of $\sB$; 
elements of $\V(P)$ are the ``ideal vertices'' of $P$. For each $\eta\in \V(P)$ and every horoball $\B\subset \H^n$ centered at $\eta$, the intersection $P\cap \B$ has finite volume. 

The convex polyhedron $P=P_\sB$ is the fundamental domain of the reflection group $\Ga_\sB< \mathrm{Isom}(\H^{n+1})$ generated by isometric reflections in facets of $P$, if and only if each dihedral angle of $P$ belongs to  $\frac{\pi}{\N}$, see \cite{Maskit, Rat}.

\medskip 
 We start the proof of Part 1 with two auxiliary lemmata.

For each nonempty closed convex subset $C\subset \H^{n+1}$ we have the {\em nearest-point projection} $\pi_C: \H^{n+1}\to C$. This projection is a continuous (actually, $1$-Lipschitz) map. 

\begin{lemma}\label{lem:L1}
Suppose that $C$ is a closed convex subset of $\H^{n+1}$ and $C\subset A$, where $A$ is a  closed  subset of $\H^{n+1}$ such that $\geo C=\geo A$. Then the restriction of $\pi_C$ to $A$ is a 
proper map, i.e. preimages of compact subsets of $C$  are compact. 
\end{lemma}
\proof  Assume that this is not the case. Then there exist sequences $y_i\in C$ and $x_i\in A$ such that $\pi(x_i)=y_i$, $\lim y_i=y\in C$, while $\lim x_i=\xi\in \geo A=\geo C$.  Let $\rho=y\xi$ denote the geodesic ray in $\H^{n+1}$ emanating from $y$ and asymptotic to $\xi$. By convexity of $C$, $\rho$ is contained in $C$. Since 
$\lim x_i=\xi$, the angles $\angle  x_i y \xi$ converge to zero as $i\to \infty$, since $x_i\to \xi$. If follows that the points $z_i:=\pi_{\rho}(x_i)$ satisfy
$$
\lim_{i\to\infty} \frac{d(x_i,z_i)}{d(x_i, y_i)}= \lim_{i\to\infty} \frac{d(x_i,z_i)}{d(x_i, y)}=0. 
$$
This is a contradiction since $z_i\in C$ and $y_i$ is the nearest-point projection of $x_i$ to $C$.  \qed 

\begin{lemma}\label{lem:L0}
Suppose that $\Gamma'\lhd \Gamma < \mathrm{Isom}(\bH^{n+1})\cong PO(n+1,1)$, where $\Gamma'$ is a Zariski dense discrete subgroup.
Then $\Gamma$ is also discrete.\footnote{The same holds for subgroups of arbitrary algebraic Lie groups with discrete center.} 
\end{lemma}
\proof Let $\gamma_i\in \Gamma$ be a sequence converging to $1\in PO(n+1,1)$. Then for every $\gamma\in \Gamma'$ we have
$$
\lim_{i\to\infty} \gamma_i \gamma \gamma_i^{-1} = \gamma. 
$$
Since $g_i:= \gamma_i \gamma \gamma_i^{-1}\in \Gamma'$ and the latter is discrete, it follows that $g_i=\gamma$ for all 
sufficiently large $i$ (depending on $\gamma$). Since $\Gamma< PO(n+1,1)$ is Zariski dense, it contains a Zariski dense Schottky subgroup $\Sigma$ of finite rank $r > 1$ with free generators $\si_1,...\si_r$. 
Hence there exists $i_0$ such that  
for all $i\ge i_0$  we have that 
$$
\gamma_i \si_k \gamma_i^{-1}=\si_k, \quad k=1,...,r. 
$$
By the Zariski density of $\Sigma$, we obtain that $\gamma_i$ belongs to the center of $PO(n+1,1)$ for all $i\ge i_0$. Since $PO(n+1,1)$ is centerless, $\gamma_i=1$  for all $i\ge i_0$. Thus, $\Gamma$ is discrete. \qed

\medskip 
Recall that $\Ga_1:=\Ga_\sB < \Isom(\H^{n+1})$ is the   discrete reflection group with the convex fundamental polyhedron $P\subset \H^{n+1}$. 
%\textcolor{red}{Maybe it makes sense to use the notation $\Ga_R$, where $R$ stands for ``reflection''} 
As we noted earlier, the limit set of $\Gamma_\sB$ is the entire sphere at infinity $\bS^{n}=\geo \H^{n+1}$;  
in particular, $\Ga_\sB$ is Zariski dense in $PO(n+1,1)$. 
%The assumption that the limit set of $\Ga_\sB$ is $\bS^{n}$ is equivalent to the hypothesis that the set of round balls bounded by the spheres in the bug $\sB$ is dense in $\bS^n$. 
Since $\sB$ is a bug, the polyhedron $P$  has infinitely many faces, equivalently,  
%We will also assume that 
$\A(P)$ is nonempty. %, equivalently, $P$ has infinite volume, equivalently, $P$.  
Thus, if we put an accumulation point in  $\A(P)$ at infinity (in the upper half-space model), then the spheres of $\sB$ are cooriented so that they bound round balls in $\R^{n}$. In this setting, 
the density condition simply means that the union of these round balls is dense in $\R^{n}$. 

The group $\Ga_2:=\Ga_S$ is a {\em symmetry group} of the bug $\sB$, i.e. a group of isometries of $\H^{n+1}$ preserving $P$. 
Then supergroup of the bug $\sB$ is the subgroup $\widetilde\Ga< \Isom(\H^{n+1})$ generated by $\Ga_1, \Ga_2$. 

Since $\Ga_2$ preserves $P$, it normalizes $\Ga_1$; hence, $\Ga_1$ is a normal subgroup of $\widetilde\Ga$.  
In view of \lemref{lem:L0} and Zariski density of $\Ga_1$, the subgroup $\widetilde\Ga< \Isom(\H^{n+1})$ is discrete. 
Our goal is to show that $\widetilde\G$ is not just Zariski dense and discrete, but is also a lattice.
We build up to the general case from some initial assumptions which simplify the exposition.

\begin{theorem}\label{thm:T1}
If, in addition to the above assumptions, $\Gamma_2$ is convex-cocompact and $\A(P)=\geo P$ (i.e. $\V(P)=\emptyset$), 
then the supergroup $\widetilde\Gamma$ is a uniform lattice in $\Isom(\H^{n+1})$. 
\end{theorem}
\proof Since $\H^{n+1}/\widetilde\Ga$ is homeomorphic (actually, isometric) to  $P/\Ga_2$,  it suffices to show that the quotient $P/\Gamma_2$ is compact.  
Since $\Ga_2$ is convex-cocompact, it acts cocompactly on the convex hull $C$ of its limit set $\La= \geo P$. Let $K\subset C$ be a compact fundamental set of the action of $\Gamma_2$  on $C$, i.e. $\Gamma_2\cdot K=C$. Then, by \lemref{lem:L1} (applied with $A=P$), the preimage $K':=\pi^{-1}_C(K)\cap P$ is also compact. Since $\pi_C$ is equivariant with respect to the action of $\Gamma_2$, it follows that $K'$ is a compact fundamental set of the action of $\Gamma_2$ on $P$. Hence, $P/\Ga_2$ is compact and  $\widetilde\Gamma$ is a uniform lattice. 
%{\misha Since $\Gamma_1\cdot P=\bH^{n+1}$, it follows that $\widetilde\Gamma \cdot K'= \bH^{n+1}$ and, hence, $\widetilde\Gamma$ is a uniform lattice.} 
\qed 

\medskip
%Below we present two generalizations of \thmref{thm:T1} in the context of (possibly nonuniform) lattices. 

Next we relax the assumption that the ideal boundary $\geo P$ coincides with the accumulation set $\A(P)$ (that is, we allow $\V(P)$ in \eqref{eq:geoP} to be nonempty), while keeping convex-cocompactness of $\G_2$.
We let $\mu_{n+1}$ denote the {\em Margulis constant} of $\H^{n+1}$ (see e.g. \cite{Kapovich-book} or \cite{Rat}).

\begin{theorem}\label{thm:T2}
If $\Gamma_2$ is convex-cocompact, then the supergroup  $\widetilde\Gamma$ is a  lattice in $\Isom(\H^{n+1})$. 
\end{theorem}
\proof Let $C=C_{\Ga_2}$ again denote the closed convex hull of the limit set of $\Ga_2$. 

Each ideal vertex $\eta_j\in \V(P), j\in J$, represents a finite volume cusp of the orbifold ${\mathcal O}= \H^{n+1}/\Gamma_1\cong P$; its stabilizer $\Ga_{1,j}$ in $\Ga_1$ is generated by reflections in the faces of $P$ asymptotic to $\eta_j$.  
We will be using certain open horoballs $\B_{\eta_j}\subset \bH^{n+1}$ centered  at the points $\eta_j$. We call such a horoball $\B\subset \bH^{n+1}$ an {\em $\eps$-Margulis horoball}, where $0< \eps \le \mu_{n+1}$, if it is the maximal horoball such that for each $x\in  \B$, there exists a parabolic element $g\in \Ga_{1,j}$ satisfying
$$
d(x, g(x))< \eps. 
$$
The Margulis lemma implies that Margulis horoballs $\B_{\eta_j}$ are pairwise disjoint. In view of \eqref{eq:inrad}, 
we can choose $\eps$ such that  each $\eps$-Margulis horoball is disjoint from $C$.\footnote{This is not really necessary but provides a clearer picture of the situation.} From now on, we fix a collection  $\B_{\eta_j}, j\in J$, of such $\eps$-Margulis horoballs. This collection is necessarily $\Ga_2$-invariant. Since each intersection $\B_{\eta_j}\cap P$ has finite volume
and $\H^{n+1}/\widetilde\Ga$ has the same volume as $P/\Gamma_2$, in order to prove that  $\widetilde\Gamma$ is a lattice it suffices to show that the quotient
$$
(P\ \setminus\  \bigcup_{j\in J} \B_{\eta_j})/\Gamma_2
$$
is compact. (This will also force finiteness of the number of $\Ga_2$-orbits of the horoballs $\B_{\eta_j}$.) 
 Compactness of the quotient again follows by  applying  \lemref{lem:L1} to the set 
$$
A:= P \ \setminus \  \bigcup_{j\in J} \B_{\eta_j} 
$$
as in the proof of \thmref{thm:T1}.  \qed 

\medskip 
Lastly, we would like to relax the convex-cocompactness assumption on $\Gamma_2$, replacing it with geometric 
finiteness. Unfortunately, in  general, the subgroup $\widetilde\Gamma< \Isom(\H^{n+1})$ will not be a lattice, already in the case $n=3$.  This can happen even if the limit set $\A(P)$ is Zariski-dense in $\bS^n$.

\begin{figure}
\includegraphics[width=.2\textwidth]{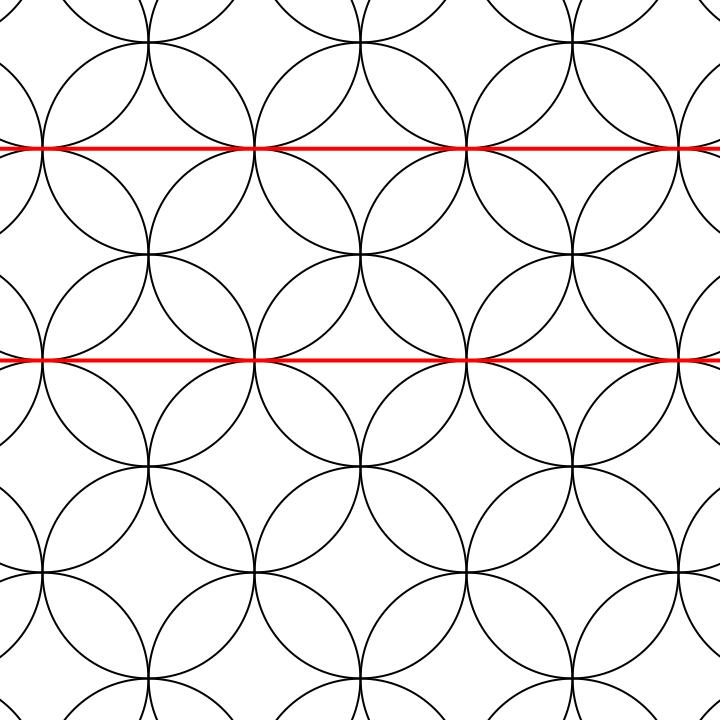}
\caption{A Kleinian bug which is not cusp-compatible.}
\label{fig:4}
\end{figure}

\begin{example}\label{ex:nonCusp}
{\em
A concrete example is the bug consisting of  circles of unit diameter centered at the set  
$
\Z[i] \cup (\foh + \foh i + \Z[i]) ,
$
that is, the Gaussian integers and their translates by the vector $\foh+\foh i$. The group $\Gamma_1$ is generated by reflections in these circles. The group $\Gamma_2$ is generated by reflections in the lines $\{\Re(z)=0\}$ and $\{\Re(z)=1\}$. The situation is illustrated in \figref{fig:4}. The accumulation point 
is $\infty$, and the group $\widetilde\G=\<\Ga_1,\Ga_2\>$ is a non-lattice. On the other hand, if we take $\Gamma'_2$ generated by the reflections in the lines 
$\{\Re(z)=0\}$, $\{\Re(z)=1\}$, $\{\Im(z)=0\}$, $\{\Im(z)=1\}$, then the group $\widetilde\G'=\<\Ga_1,\Ga'_2\>$ is a lattice. 
}
\end{example}

\begin{example}
{\em
More generally,
let $\Gamma_2$ be a geometrically finite subgroup of $\Isom(\H^{3})$ which has at least one rank-one cusp. Then (by  either using the Andreev--Thurston theorem, see \cite[\S\S 13.6--13.7]{Thurston},  or Brooks's theorem, see \cite{Brooks86} or \cite{Kapovich-book}) there exists a quasiconformal deformation $\Gamma'_2$ of $\Gamma_2$ and a reflection group $\Gamma_1$ such that the fundamental domain of $\Gamma_1$ is invariant under $\Gamma'_2$ and its accumulation set equals the limit set of $\Gamma'_2$, but the groups $\Gamma_1$, $\Gamma'_2$ are not cusp-compatible.  
}
\end{example}

A necessary condition for $\widetilde\Ga$ to be a lattice is 
that the $\widetilde\Ga$-stabilizer of each parabolic fixed point of $\Ga_2$ is virtually abelian of  virtual rank $n$. 
(Recall that a group $G$ is a virtually 
%free 
abelian if it contains a free abelian subgroup 
%$A'$ 
of finite index. The virtual rank 
of $G$ is defined to be the 
rank of any 
maximal
free abelian subgroup of $G$. By Bieberbach's Theorem, a discrete parabolic subgroup of $\Isom(\H^{n+1})$ is necessarily virtually abelian of virtual rank $\le n$.)  We thus arrive at the following definition of cusp-compatibility, for 
the pair of discrete groups $\Gamma_1=\Gamma_{\sP}$, $\Gamma_2=\Gamma_S$, repeated from \defref{def:cuspComp}.

\begin{definition}
The subgroups $\Ga_1, \Ga_2< \Isom(\H^{n+1})$ are {\em cusp-compati\-ble} if for every maximal parabolic subgroup $\Pi_2< \Gamma_2$, there exists a subgroup  $\Pi_1< \Gamma_1$  such that the subgroup $\Pi$ generated by $\Pi_1, \Pi_2$ is virtually abelian of virtual rank $n$. 
\end{definition}

\medskip 
Observe that the subgroup $\Pi_1$ in this definition necessarily fixes the limit set (the single fixed point at infinity) $\{\la\}\subset\bS^n$ 
of $\Pi_2$. Thus, $\Pi$ fixes $\la$ as well. In view of the discreteness of $\widetilde\Gamma$ (and, hence, of $\Pi$), we can (and will) as well assume that $\Pi_1$ is the full stabilizer of  $\la$ in $\Ga_1$. (The group $\Pi_1$ might be finite.) Since $\Ga_1$ is normal in $\widetilde\Ga$, the subgroup $\Pi_1$ is normal in $\Pi$ and $\Pi$ splits as the semidirect product $\Pi_1\rtimes \Pi_2$.

\begin{lemma}\label{lem:domains}
1. The group $\Pi_1$ is generated by reflections in the facets of $P$ asymptotic to $\la$ and, hence,  
has a fundamental domain $P_\la$ in $\H^{n+1}$ equal to the intersection of all 
half-spaces defined by facets of $P$ asymptotic to $\la$. 

3. The fundamental domain $P_\la$ is invariant under the action of $\Pi_2$. 
\end{lemma}
\proof 1. 
Let $s_\alpha$ denote the isometric reflections in the facets $F_\alpha$ of $P$ generating the reflection group $\Gamma_1$. 
Every element $\gamma\in \Gamma_1$ is represented by a reduced word $w$ in the generators $s_\alpha$. 
Suppose that $w$ has the form $u \cdot s_\alpha \cdot v$, where $v$ is the product of reflections in the facets of $P_\la$, while $s_\alpha$ is the reflection in a facet  $F_\alpha$ not asymptototic to $\lambda$. Then $s_\alpha v(P)$  is contained in the half-space $H^-_\alpha$ whose closure in $\H^{n+1}\cup \geo \H^{n+1}$ does not contain $\la$. It follows that 
$u\cdot s_\alpha \cdot v(P)$ is also contained in $H^-_\alpha$ and, hence, $\gamma$ cannot possibly fix $\la$. 

2. Since $P$ is preserved by $\Ga_2$, and, hence, by its subgroup $\Pi_2$, and $\Pi_2$ fixes $\la$, the elements of $\Pi_2$ send facets of $P$ asymptotic to $\la$ to facets of $P$ asymptotic to $\la$. \qed  

\medskip 
The following %(straightforward) 
lemma provides a list of equivalent algebraic and geometric characterizations of 
 cusp-compatibility in the context of the pair of groups $\Gamma_1=\Gamma_{\sP}$, $\Gamma_2=\Gamma_S$: 

\begin{lemma}\label{lem:L2}
The following are equivalent for  subgroups $\Pi_1, \Pi_2$ of $\Ga_1, \Ga_2$ as above: 

\begin{enumerate}
\item The subgroup $\Pi$ generated by $\Pi_1, \Pi_2$ is virtually abelian of virtual rank $n$. %1

\item Virtual ranks of $\Pi_1, \Pi_2$ add up to $n$. %2

\item $\Pi$ acts cocompactly on $\geo \H^{n+1} \ \setminus\ \{\la\}$. %3

\item  $\Pi_2$ acts cocompactly on the intersection $P_\la\cap \partial \B$, for every  horoball $\B\subset \H^{n+1}$ centered at $\la$. %4

\item  $\Pi_2$ acts cocompactly on the intersection $P\cap \partial \B$, for every  horoball $\B\subset \H^{n+1}$ centered at $\la$. %5

\item  $\Pi_2$ acts with finite covolume on the intersection $P_\la\cap \B$, for every  horoball $\B\subset \H^{n+1}$ centered at $\la$. %6

\item $\Pi_2$ acts with finite covolume on the intersection $P\cap \B$, for some  horoball $\B\subset \H^{n+1}$ centered at $\la$. %7
\end{enumerate}
\end{lemma}
\proof The equivalence (1) $\Leftrightarrow$ (2) follows from the semidirect product decomposition $\Pi=\Pi_1\rtimes \Pi_2$. 
The equivalence (1) $\Leftrightarrow$ (3) follows from the fact that a discrete subgroup of $\R^n$ is a uniform lattice  
 if and only if this subgroup has rank $n$.  

To prove the equivalence (3) $\Leftrightarrow$ (4) observe that, by  \lemref{lem:domains},  
$P_\la\cap \partial \B$ is the fundamental domain for the action of $\Pi_1$ on the horosphere $\partial \B$ and this intersection is 
invariant under $\Pi_2$. Therefore, cocompactness of the action of $\Pi$ on the horosphere is equivalent to the cocompactness 
of the action of $\Pi_2$ on $P_\la\cap \partial \B$. 

The proofs of equivalences  (3) $\Leftrightarrow$  (5) $\Leftrightarrow$  (6)$\Leftrightarrow$  (7) are similar to the proof of 
 (3) $\Leftrightarrow$ (4)  and are left to the reader.  \qed

\medskip 
Subgroups $\Pi_1, \Pi_2$ are called {\bf compatible} if they satisfy one of the equivalent conditions in this lemma.   
Note that since $P\subset P_\la$,  compactness of $(P_\la\cap \partial \B)/\Pi_2$, implies that only finitely many $\Pi_2$-orbits 
of faces of $P$ might intersect $\B$. Moreover, the compatibility of $\Pi_1, \Pi_2$ (in the form of the 3rd condition) 
implies that any sequence of faces of $P$ converging to $\la$ (possibly outside of a horoball $\B$) 
is contained in finitely many $\Pi_2$-orbits of faces. 

 \medskip
 
 \begin{lemma}\label{lem:commensurable}
 Suppose that $\A(P)$ has cardinality $\ge 2$. Then any two symmetry groups of $\sB$ are commensurable and, in particular, 
 the  cusp-compatibility of the subgroups $\Ga_1, \Ga_2$ depends only on $P$.  
 \end{lemma}
\proof Suppose that $\Ga_2, \Ga'_2< \Isom(\H^{n+1})$ are geometrically finite subgroups preserving $P$ and having $\A(P)$ as their limit sets. 
Since  $|\A(P)|\ge 2$, by the definition of geometric finiteness, the groups $\Ga_2, \Ga'_2$ act with finite covolume on the $\eps$-neighborhood $Y$ of the convex hull of 
%the limit 
$\A(P)$ (which is their common limit set). 
Since both $\Ga_2, \Ga'_2$ preserve the polyhedron $P$, they generate a discrete subgroup 
$\Ga_3< \Isom(\H^{n+1})$, which, therefore, acts properly discontinuously, with finite covolume on $Y$. It follows that $|\Ga_3: \Ga_2|<\infty$,  
$|\Ga_3: \Ga'_2|<\infty$, which implies commensurability of $\Ga_2, \Ga'_2$. Therefore, the pair $(\Ga_1, \Ga_2)$ is cusp-compatible if and only if $(\Ga_1, \Ga'_2)$ is. \qed 

\begin{rem}{\em
 \exref{ex:nonCusp} shows that, when  $|\A(P)|=1$, cusp-compatibi\-lity depends not only on $P$ 
but also on $\Ga_2$. 
}
\end{rem}

\begin{thm}\label{thm:T3}
If $\Ga_2$ is geometrically finite and the pair $\Gamma_1, \Gamma_2$ is cusp-compatible, 
then the subgroup $\widetilde\Gamma$ generated by $\Gamma_1, \Gamma_2$  is a lattice in $\Isom(\H^{n+1})$. 
\end{thm}
\proof Similarly to the proof of \thmref{thm:T2}, we define open horoballs $\B_{\eta_j}, j\in J$. Since 
each $\eta_j$ does not belong to $\A(P)$, it follows that all $\eta_j$'s lie in the discontinuity domain of 
the group $\Ga_2$. This, together with geometric finiteness of $\Ga_2$ implies that 
there exists a collection $\{\bar{\B}_{\la_i}; i\in I\}$ of 
 closed (Margulis) horoballs centered at parabolic fixed points of $\Ga_2$, such that:

\begin{enumerate}
\item The collection of horoballs is $\Ga_2$-invariant. 

\item The horoballs are pairwise disjoint and disjoint from the horoballs $\bar{\B}_{\eta_j}$. 

\item Each horoball $\bar{\B}_{\la_i}$ intersects only those faces of $P$ which are asymptotic  to $\la_i$. 
\end{enumerate} 

Set
$$
A:= P \ \setminus\ \left( \bigcup_{j\in J} \B_{\eta_j} \cup \bigcup_{i\in I} \B_{\la_i}\right).  
$$

There are two cases to consider. First, suppose that $\A(P)=\La(\Ga_2)$ is a singleton $\{\la\}$. 
Then $\{\B_{\la_i}: i\in I\}$ consists of a  single horoball $\B$. Consider the nearest-point projection $\pi_{\bar{\B}}: \H^{n+1}\to \bar{\B}$. 
For $x\notin \bar{\B}$, this projection is given by the unique intersection point of the geodesic ray $x\la$ with the horosphere $\partial \B$. As in the proof of \thmref{thm:T2}, the restriction of $\pi_{\bar{\B}}$ 
to $A$ is a proper map. The image $\pi_{\bar{\B}}(A)$ is contained in $\partial \B\cap P$: By the convexity of $P$, for each 
$x\in A$ the ray $x\la$ is contained in $P$, hence, the intersection point $x\la\cap \partial \B$ belongs to 
$\partial \B\cap P$. 

Since  the pair $\Gamma_1, \Gamma_2$ is cusp-compatible, by \lemref{lem:L2}(5), the action of 
$\Ga_2$ on $P\cap \partial \B$ is cocompact. Combining this with the fact that the map $\pi_{\bar{\B}}: A\to \bar{\B}$ is proper and its image 
is contained in $\partial \B\cap P$, we conclude that the quotient  $A/\Gamma_2$ is compact. 
Since $(P\cap \B)/\Ga_2$ has finite volume and the volume of each cusp $P\cap \B_{\eta_j}$ is finite, 
we conclude that $P/\Ga_2$ has finite volume and, therefore, $\widetilde\Gamma$ is a lattice. 

\medskip 
We now consider the ``generic'' case, when $\Gamma_2$ has at least two limit points. As before, let $C=C_{\Ga_2}$ denote 
the closed convex hull  in $\H^{n+1}$ of the limit set $\La=\La(\Ga_2)$. Since $\La= \A(P)\subset\geo P$, 
convexity of $P$ implies that $C\subset P$. We let $C'\subset C$ be the complement to the union of open horoballs $\B_{\la_i}, i\in I$. The group $\Gamma_2$ acts cocompactly on $C'$ since it is geometrically finite. 
As in the proof of \thmref{thm:T2}, the restriction of $\pi_C$ to 
$A$ is a proper map and $(A\cap \pi_C^{-1}(C'))/\Ga_2$ is compact. Lastly, \lemref{lem:L2} (and the discussion following it) implies that for every horoball $\B_{\la_i}, i\in I$, its $\Ga_2$-stabilizer $\Pi_{\la_i}$  acts with finite covolume on 
$$
P\cap \pi_C^{-1}(\B_{\la_i}).  
$$
Hence, $P/\Ga_2$ has finite volume and, since $Vol(P/\Ga_2)= Vol(\H^{n+1}/\widetilde\Ga)$, it follows that 
$\widetilde\Ga$ is a lattice. \qed 

\medskip
The next lemma establishes the cusp-compatibility condition in the 
the special case of bugs which are packings. % instead of more general bugs.
%%%%
%\begin{comment}
% ``simple'' bugs analogously to \cite{KN}. 
%\begin{definition}\label{def:simple}
%A Kleinian bug $\sB$ is called {\em simple} if the hyperplanes defining facets of $P$ are pairwise disjoint, equivalently, the open balls defined by the cooriented spheres in the bug are pairwise disjoint.
%\end{definition}
% That is, a bug is simple if it is a sphere packing, but with the ``Kleinian'' condition replacing ``crystallographic.''
%\end{comment}
%%%%%

\begin{lemma}\label{lem:simple->lattice}
If $\sB$ is a Kleinian packing then the cusp-compatibility condition always holds. 
\end{lemma}
\proof Since $\sB$ is a Kleinian packing, the polyhedron $P$ is the intersection of half-spaces $H^+_\alpha$, which are 
bounded by pairwise disjoint hyperbolic hyperplanes $H_\alpha$. In particular, $\geo P$ has no isolated points 
and, hence, equals the accumulation set $\A(P)$, which is the limit set $\La$ of $\Ga_S$, and $P=C_{\Ga_S}$, the closed 
convex hull of $\La$. Since $\Ga_S$ is geometrically finite, $C_{\Ga_S}/\Ga_S$ has finite volume. But  
$C_{\Ga_S}/\Ga_S=P/\Ga_S$ is naturally isometric to $\H^{n+1}/\widetilde\Ga$. We conclude that $\widetilde\Ga$ is a lattice, which is 
equivalent to the cusp-compatibility condition. \qed 

\

{\bf Part 2.}\

Proof of Part 2(i). Step 1. For each  defining hyperplane $H_\al$ of $D$, $\al\notin A$, we let $H_\al^+\subset \H^{n+1}$ denote 
 the closed half-space  bounded by $H_\al$  and containing $D$. We define $D_S$ as the intersection of these half-spaces, 
$$
D_S=\bigcap_{\al\notin A} H_\al^+.
$$
Each facet of $D_S$ is contained in one of the facets of $D$ which is not of the form $F_\al, \al\in A$. And conversely, each facet of $D$ not of the form 
$F_\al, \al\in A$, is contained in one of the facets of $D_S$. The generators $\gamma\in S$ still pair the facets of $D_S$. We leave it to the reader to check that, since $D$ is a fundamental polyhedron of $\Gamma$, the conditions of the Poincar\'e's Fundamental Domain Theorem (see \cite{Rat}) still hold for $D_S$ and the face-pairing transformations $\gamma\in S$, and, hence, $D_S$ is a fundamental polyhedron of $\Gamma_S$. Since, by the construction, $D_S$ is finitely-sided, the group $\Gamma_S$ is geometrically finite with the fundamental domain $D_S$.

\medskip 
Step 2. We let ${\mathcal R}$ denote the collection of all $\widetilde\Gamma$-images of the hyperplanes $H_\al, \al\in A$. Since $D$ is a fundamental domain of $\widetilde\Ga$, none of the hyperplanes in ${\mathcal R}$ intersects the interior of $D$. The hyperplanes in ${\mathcal R}$ define a partition of $\H^{n+1}$ into (convex) connected components; the closure of one of these components, denoted $P$, contains $D$.  We let $\Gamma_R< \widetilde\Gamma$ denote the subgroup generated by reflections in the hyperplanes $H\in {\mathcal R}$. Thus, 
$\Gamma_R$ is a nontrivial normal subgroup in the lattice $\widetilde\Ga$; hence, the limit set of $\Ga_R$ is the entire sphere $\bS^n$. The polyhedron $P$ 
 is a fundamental domain of $\Gamma_R$ and  $\Gamma_R$ is generated by reflections in those hyperplanes $H\in {\mathcal R}$ for which $H\cap P$ are 
facets of $P$. Since $\Gamma_R< \widetilde\Ga$ and the latter is a lattice, it follows that the dihedral angles of $P$ come from a finite subset of $\frac\pi\N$. 
We conclude that $P$ defines a bug $\sB$. 

\medskip 

Step 3. We next verify that $\Ga_S$ is a symmetry group of the bug $\sB$. We first check that the generators $\gamma\in S$ of $\Ga_S$ preserve the polyhedron $P$. It suffices to show that each $\gamma$ sends facets of $P$ to facets of $P$. The element $\gamma$ pairs a facet $F_\gamma$ of $D$ to a facet $F'_\gamma$ of $D$. 
Let $G_\beta$ be a facet of $P$ and $H_\beta\in {\mathcal R}$ be the hyperbolic hyperplane containing $G_\beta$. 
We pick a generic base-point $o$ in the interior of $D$ and  a generic point $x$ on the facet $G_\beta$. Then, by the convexity of $P$ and since $o\in D\subset P$, the geodesic segment $ox$ is disjoint from all the hyperplanes $H\in {\mathcal R}$ except for the point $x\in G_\beta$. Similarly, the geodesic segment $o \gamma(o)$ is contained in $D\cup \gamma(D)$ and crosses their intersection at an interior point of $F'_\gamma$ (since $o$ was chosen generically). It follows that the 
segment $o \gamma(o)$ is also disjoint from all the hyperplanes in ${\mathcal R}$. Thus, the union
$$
o \gamma(o) \cup \gamma(o) \gamma(x)
$$
is a path connecting $o$ to $\gamma(x)$ and disjoint from all the hyperplanes in ${\mathcal R}$ except for the point $\gamma(x)\in \gamma H_\beta$. Hence, $\gamma(x)$ lies in a facet of $P$. Thus, we verified that the generators $\gamma\in S$ have the property that they send facets of $P$ to facets of $P$ and, moreover, respect their coorientation: 
The half-space $H^+_\al$ determined by a facet $F_\alpha$ and containing $o$, maps to the half-space $H^+_{\al'}$ containing $o$. Therefore, the entire group $\Gamma_S$ preserves $P$. 

\medskip 

Step 4. Since, by the construction, $\widetilde\Ga$ is a lattice generated by $\Ga_R$ and $\Ga_S$, the group $\Gamma_S$ acts on the set of facets of $P$ with finitely many orbits. 
Thus, if 
$\xi$ is an accumulation point of the bug $\sB$, there is a facet $G_\al$ of $P$ and an infinite sequence $\gamma_i\in \Gamma_S$ such that $\gamma_i(G_\al)$ converges to $\xi$. In other words, ${\mathcal A}(\sB)$ is contained in the limit set of $\Gamma_S$. The opposite inclusion follows from the fact that $\Gamma_S$ preserves the polyhedron $P$. We conclude, therefore, that $\Gamma_S$ is a symmetry group of $\sB$ and $\sB$ is geometric with $\Gamma_R=\Gamma_{\sB}$. 

\medskip 

Step 5. Lastly, the cusp-compatibility of the pair of groups $\Gamma_R, \Gamma_S$ follows from the fact that they generate a lattice $\widetilde\Gamma$. 

Thus, we proved that $\sB$ is a Kleinian bug, which concludes the proof of Part 2(i) of  \thmref{thm:lattice}. 

\medskip 
Proof of Part 2(ii). We continue with the notation from the proof of Part 2(i). 
However, before starting the actual proof of (ii) we will  have a discussion related to combinatorics of 
convex  fundamental polyhedra. Each facet $F$ of $\partial D$ is {\em paired} with 
another facet $F'$ by a unique generator $\ga=\ga_{F,F'}\in S'$; $\ga(F)=F'$, where, possibly, $F=F'$.   

\begin{rem}{\em
Following Ratcliffe in \cite{Rat}, we require that if a generator $\ga\in S'$ preserves a facet $F$ of $D$, then it fixes $F$ pointwise, i.e. is a reflection in $F$. To achieve this, one performs, if necessary, a subdivision of geometric facets of $D$. We refer to 
\cite{Rat} for details. Accordingly, if $F_1, F_2$ are distinct facets of the same geometric facet of $\partial D$ and 
$\ga_i: F_i\to F'_i, i=1,2$, pair these facets, and it happens that $\ga_1=\ga_2$ are equal as elements of $\widetilde\Ga$, then we still regard $\ga_1, \ga_2$  as distinct elements of $S'$.   
}
 \end{rem}

 A {\em ridge} of a convex polyhedron in $\H^{n+1}$ is the $n-1$-dimensional intersection of two facets. 
 %A pair $(F, E)$, where $F$ is a facet of $D$ and $E$ is its boundary ridge, is a {\em flag}.  
 
 \medskip
 {\bf The pseudogroup ${\mathcal G}$.}  The pair $(D, S')$ defines a pseudogroup ${\mathcal G}$  acting on $\partial D$, which we discuss below. Each generator $\gamma\in S'$ is an element of 
 ${\mathcal G}$; it is a partially defined map between the two facets of $D$ paired as $\ga: F\to F'$. 
 Then the unique generator sending $F'\to F$ is $\ga^{-1}$.  A composition $\ga_2\circ \ga_1$ of two generators 
 $$
  \ga_1: F_1\to F_1', \ga_2: F_2\to F_2'
 $$
 is {\em admissible}  if $F_1'\cap F_2$ is a common ridge of these facets. A (possibly empty) word 
 $$
w= \ga_l\circ ...\circ \ga_1, \ga_i: F_i\to F'_i, i=1,...,l,  
 $$
 is {\em admissible} if each consecutive composition in it is admissible. Thus, each admissible word defines a map from one ridge $E$ 
 to another ridge $E'$, the domain and the range of $w$ (unless $l=1$ in which case the domain and the range are facets). 
 Here $E$ is a boundary ridge of $F_1$ and $E'$ is a boundary ridge of $E_l'$.  The pseudogroup ${\mathcal G}$ then consists of admissible compositions of the generators. Note that each admissible word is necessarily a reduced word in the alphabet $S'$: 
 $\ga_{i+1}\ne \ga^{-1}_i$ for each $i$. 
 
The pseudogroup ${\mathcal G}$ defines an equivalence relation $\sim_{{\mathcal G}}$ on $D$:
$x\sim_{{\mathcal G}} y$ ($x$ is ${\mathcal G}$-equivalent to $y$) if and only if there exists an element $\ga\in {\mathcal G}$ sending $x$ to $y$.  This equivalence relation is the one obtained by saturating the non-reflexive and non-transitive relation given by
$$
x\sim y, x\in F, y\in F', y=\ga_{F,F'}(x).
$$
 An important fact, coming from the assumption that $D$ is a fundamental polyhedron of $\Ga$ is that the natural projection of  
 quotient spaces 
 $$
 D/_{\sim_{\mathcal G}} ~ \to {\mathbb H}^{n+1}/\widetilde\Ga$$ is a homeomorphism; 
 this implies that two points in $D$ are ${\mathcal G}$-equivalent if and only if they are $\widetilde\Ga$-equivalent, i.e. belong to the same  
 $\widetilde\Ga$-orbit. 
 
 \medskip 
 {\bf Ridge-chains and cycles.} Suppose that $E=E_1= F_1\cap F_0'$ is a ridge of $D$. Let $\ga_1\in S'$ denote the generator  pairing the facet $F_1$ to another 
 facet, $F_1'$; this yields a new ridge $E_2:=\ga_1(E_1)=F_1'\cap F_2$. Then let $\ga_2\in S'$ be the generator pairing $F_2$ to a facet $F_2'$. 
 The composition $\ga_2\circ \ga_1$ is admissible (by the construction). 
  This composition process continues (uniquely),  until we return to the original ridge so that  
 $$
 \ga_k: E_k\to E_1, \quad \hbox{and}\quad \ga_k(F_k)=F'_k= F_0'. 
 $$
 (The process has to terminate since $D$ has only finitely many faces; in the case of polyhedra with infinitely many faces, such termination is a consequence of one of the axioms of fundamental polyhedra.) There is an important caveat regarding 
 this definition that applies in the special case when the {\em first return} to the initial ridge yields $F'_k=F_1$ instead of $F'_k=F'_0$;  we discuss this in  \rmkref{rem:repetition} below. 
  
 The finite sequence
 $$
c_E= (\ga_1, \ga_2,...,\ga_k)
 $$
 is called a {\em ridge-cycle}; it corresponds to the word 
 $$
w_{c_{E}}=\ga_k\circ ...\circ \ga_1,
$$
which is an admissible composition. Its subwords 
  $$
w_{c_{E,E'}}=\ga_l\circ ...\circ \ga_1,\quad l\le k,
$$
correspond to {\em ridge-chains}
$$
c_{E,E'}= (\ga_1, \ga_2,...,\ga_l),
$$ 
where $E'=E_{l+1}=\ga_l(E_l)$. (The notation is slightly ambiguous since  the chain 
$c_{E,E'}$ is  not uniquely determined by $E, E'$; the same, of course, applies to the notation $c_E$.) 
The  element $\ga_{c_{E,E'}}\in \Ga$ corresponding to the composition $w_{c_{E,E'}}$ sends $E_1$ to $E_{k+1}$.  
We will refer to  $w_{c_{E,E'}}$ as the {\em word of the ridge-chain} $c_{E,E'}$.  
 
 We let $\theta_i$ denote the interior dihedral angles of $D$ along the ridges $E_i$; the sum 
$$
\theta_{c_{E,E'}}=\sum_{i=1}^l \theta_i
$$ 
is the {\em total angle} of the ridge-chain $c_{E,E'}$.

 \begin{rem}\label{rem:repetition} 
 {\em
1. In the sequence of facets given by a ridge-chain, we could have $F_i'=F_i$; this happens when $\ga_i$ 
 is the reflection in the facet $F_i$. Accordingly, in the case $F_i'=F_i$ we will have $\ga_i(E_i)=E_{i+1}=E_i$. 
 
 2. In this situation (i.e. $F_i'=F_i$), the word $w_{c_E}$ contains a prefix subword which is a ``palindrome'' 
 $$
u= (\ga_{i-1}\circ ...\circ \ga_1)^{-1} \circ \ga _i \circ (\ga_{i-1}\circ ...\circ \ga_1). 
 $$ 
 The word $u$ represents an element of $\Ga$ sending $F_1$ back to itself and the ridge $E$ back to itself.  In this case, of course, $u$ does not send $F_1$ to $F_0'$ as required by the definition of a ridge-cycle. Thus, the actual cycle-word $w$ will be longer than $u$ and will 
 equal $v\circ u$, where $v$ is another  palindromic composition, starting with a face-pairing  $F'_0\to F_0$. 
 As a simple example of this situation, one can take the case when both $F_0', F_1$ are reflective facets with the corresponding generating reflections $\ga_0, \ga_1$ respectively. Then $w=\ga_0\circ \ga_1$. 
  
  3. Each ridge $E$ of $D$ defines exactly two ridge-cycles, which differ by swapping the facets $F_1, F_0'$, reversing the order in the sequence $(\ga_i)$ and inverting the generators in the cycle.   
 The total angle is, of course, independent of  which of the two ridge-cycles is used. 
}
  \end{rem}
 
 Conversely, suppose we are given an admissible composition 
 $$
 w=\ga_l\circ ...\circ \ga_1
 $$ 
 representing $\ga\in\widetilde\Ga$ 
 sending a ridge $E=F_0'\cap F_1$ to a ridge $E'$.  This alone, however, does not guarantee that corresponding sequence of generators is a chain since for some $i$, $1\le i < k$, we may have that $\ga_i\circ ...\circ \ga_1$ sends $E$ to itself and $\ga_i: F_i\to F'_i=F'_0$. Taking the minimal $i$ with this property we obtain a ridge-cycle $c_E$ 
 and the corresponding word $w_{c_E}$. Hence, we decompose $w$ as 
 $$
w = w'\circ (w_{c_E})^t,
 $$
  where $w'$ is a subword in $w_{c_E}$ which  represents an element $\ga'$, $\ga'(E)=E'$. In particular, $w'$ is the word of 
  a ridge-chain $c_{E,E'}$.

\medskip 
The fact that $D$ is a convex fundamental polyhedron of $\widetilde\Ga$ implies that for every ridge $E$, 
the total angle $\theta_E$ is of the form $\frac{2\pi}{m}$, 
 and $\ga_{c_E}$ has order $m$. Suppose now that the ridge-cycle $c_E$ is such that the consecutive facets 
$F_{i}, F'_{i}$ are equal (this happens when $\gamma_i$ is a reflection). In this case   
the angle $\theta_i$ appears twice in the sum defining the total angle $\theta_E$.  In particular, if 
$\theta_1=\theta_k=\frac{\pi}{2}$ and $\theta_s=\frac{\pi}{2}$ for some $1\le s\le k$, and $\gamma_s, \gamma_k$ are reflections, 
then both $\theta_1, \theta_s$ contribute {\em twice} to the total angle and, thus,   there are exactly two possibilities for 
the ridge-cycle $c_E$:

1. $k=2, s=1$, $E=E_1=E_2$, $\theta_{c_E}=\pi$,  and the ridge-cycle is as in the  simple 
example in  \rmkref{rem:repetition}(2) above. In particular, the reflections $\gamma_1, \gamma_2$ commute and 
$\gamma_i(H_j)=H_j$, $j=1, 2$ is taken modulo $2$. 

2. $k=4, s=2$, $E_1\ne E_2$, $\theta_{c_E}=2\pi$, and, up to inversion, 
$$
w_{c_E}=\tau_0\circ \gamma_1^{-1}\circ \tau_2\circ  \gamma_1$$ 
where $\tau_2=\ga_2$ is the reflection in the facet $F_2$, while $\tau_0=\ga_4$ is the reflection in the facet $F_0'=F_4$.
The composition $\ga_{c_E}$ then has order $1$, i.e. the word $w_{c_E}$ represents the neutral element of the group $\widetilde\Ga$. 
In particular, in this case, 
$$
\tau_2= \gamma_1\circ \tau_0\circ  \gamma_1^{-1}.
$$
%Thus, $\gamma_1(H_2)=H_2$ where $H_i$ is the hyperbolic hyperplane fixed by $\tau_i$, $i=1,2$. 

 In both cases 1 and 2, there are exactly four images of $D$ (one of which is $D$ itself) under the elements of $\widetilde\Ga$,  sharing the ridge $E$,  all with  right dihedral angles at $E$. The hyperplanes $H_0, H_1$ bounding these images and containing the faces $F'_0, F_1$ are orthogonal to each other and divide the hyperbolic space in four quadrants, each containing one of the above images of $D$. The stabilizer of $E$ in $\widetilde\Ga$ preserves both hyperplanes.

\medskip 
We are now ready for the proof of Part 2(ii). Suppose that the bug $\sB$ is not a packing and two walls in ${\mathcal R}$ have nonempty transversal intersection in the hyperbolic space; the intersection necessarily has codimension $2$. Since $D$ is 
a fundamental polyhedron of $\widetilde\Ga$, this intersection comes from $D$ in the following sense:  
There exists a pair of ridges $E, E'$ (possibly equal!) formed, respectively, by pairs of facets $F_1, F'_0$ (the ridge $E$) and 
$F_s, F'_s$ (the ridge $E'$). The facets $F'_0, F'_s$ are contained in hyperplanes $H=H_0, H'=H_s\in {\mathcal R}$ 
fixed by reflections $\tau_0, \tau_s$. There exists an element $g\in \widetilde\Ga$ which carries $E$ to $E'$ and sends 
$H$ to a hyperplane $g(H)$ which meets $H'$ orthogonally along $E'$.  In particular, $E, E'$ are in the same  
$\widetilde\Ga$-orbit. It follows that there exists a chain $c_{E,E'}$ corresponding to a word $w_{E,E'}$ representing an element 
$\ga=\ga_{E,E'}$ which sends $E$ to $E'$. The element $\ga$ need not be equal to $g$, but $g=\gamma'\circ \gamma$, where 
$\gamma'\in\widetilde\Ga$ is an elliptic isometry preserving $E'$. Since $H, H'$ are both in ${\mathcal R}$, the dihedral angles of $D$ along $E, E'$ are both right angles. Thus, the discussion above regarding the ridge-chain $w_{E,E'}$  applies. In particular, 
%$\gamma'$ preserves both hyperplanes $H, g(H)$ and 
both $H, g(H)$ define facets of $D$. 
Recall that by the assumptions of Part 2(ii),  the hyperplanes in ${\mathcal R}$ defining facets of 
$D$ are pairwise disjoint, which implies that $g(H)=H'$, a contradiction. %One can streamline this argument a bit... 

This concludes the proof of  \thmref{thm:lattice}.  \qed

%\newpage

\section{(Sub)Arithmeticity}\label{sec:arith}

In this section, we %make crucial use of inversive coordinate systems to 
prove the various (sub)arithmeticity theorems. We will be repeatedly using the {\em Lorentzian model} of  hyperbolic space $\bH^{n+1}$ 
and the corresponding parameterization of round spheres in $\bS^n=\geo\bH^{n+1}$ by unit vectors of the associated quadratic form.

\begin{definition}\label{defn:hyperbolic-form} 
A quadratic form $\cF$ over a totally real number field $k$ is called {\em hyperbolic} if it has signature $(n+1,1)$ in the identity embedding $k\to \R$, and is definite in all others embeddings.  
A {\em quadratic space} is an $(n+2)$-dimensional real vector space $V$ together with a real quadratic form $\cF$ on $V$ defining a bilinear from $\<\cdot,\cdot\>$ in the usual way. 
A quadratic space is said to be {\em hyperbolic} if $\cF$ has signature $(n+1,1)$.
\end{definition}

We let  $Q$ denote the {\em standard} hyperbolic quadratic form with half-Hessian: 
\be\label{eq:Qdef} 
Q \ = \
\bp
&&-\foh\\
&I&\\
-\foh&&
\ep
.
\ee

\subsection{Inversive coordinates} \label{sec:inversive} \

Before embarking on the proofs,
we recall the very convenient (in this context) {\it inversive coordinate system} (see, e.g., \cite{KtoD, LagariasMallowsWilks2002, Wilker1982}). 
For convenience, we work here with the standard hyperbolic quadratic form $Q$, but the entire discussion applies to all quadratic forms after an appropriate change of coordinates.

To a cooriented round sphere $S$ in the boundary $\bS^n=\geo\bH^{n+1}= \R^n \cup \{\infty\}$ %of the upper half space model 
having center $z=(x_1,\dots,x_n)$ and signed radius $r$, we associate the column vector
$$
v_S \ := \ (\tfrac1r,\tfrac zr,\tfrac1{\hat r})^t \ \in \ \R^{n+2}.
$$
Here $\hat r$ is the co-radius, defined to be the signed radius of the image of $S$ under reflection through the unit sphere; more concretely, 
\begin{equation}\label{eq:coradius}
\hat r={r\over |z|^2-r^2}.
\end{equation}
When $S$ is a hyperplane, the inversive coordinates are those obtained from a limit of spheres. That is, $1/r=0$, $1/\hat r$ is half the distance from the hyperplane to the origin, and $z/r$ is the unit normal to the cooriented hyperplane.

Rewriting \eqref{eq:coradius}  as
$$
|\tfrac zr|^2 - \tfrac1r\tfrac1{\hat r}=1,
$$
we see that $Q(v_S)=1$.

The M\"obius group $Mob_n$ (the group of M\"obius transformations of $\bS^n$) acts on the space $Sph(n)$ of (cooriented) round spheres and  
the group $O_Q$ of automorphisms of the form $Q$ acts on the 1-sheeted hyperboloid $\{Q(v)=1\}$.  

\begin{lemma}\label{lem:jMap}
The map $j: S\mapsto v_S$ is equivariant with respect to the actions of $Mob_n$ and $O_Q$. More precisely, $j$ conjugates the M\"obius action on the space of cooriented 
spheres $S\in Sph(n)$ to the Lorentzian action on the inversive coordinate vectors $v_S\in \{ Q_S=1\}$.  
\end{lemma}
\proof The proof is essentially contained in the proof of \cite[Theorem 7.5]{Iversen}. Iversen constructs an equivariant map $\iota$ of $\bS^n$ to the 
projectivization of the conic $\{Q=0\}$. In the proof he verifies that 
$$
\iota R_S \iota^{-1} = \tau_{j(S)},
$$
where $R_S$ is inversion in $S\in Sph(n)$ and $\tau_v$ is the Lorentzian reflection in $v^\perp$, the Lorentzian orthogonal complement to the vector $v$. 
To prove equivariance of  $j$ it remains to observe that for every $g\in Mob_n$, $S\in Sph(n)$ and $v_S=j(S)$, 
$$
\iota \sigma_{gS} \iota^{-1}= \iota g R_S g^{-1} \iota^{-1}= g \iota R_S \iota^{-1} g^{-1}= g \tau_{v_S} g^{-1}= \tau_{g v_S}.
$$ 
Hence, $j$ sends the sphere $gS$ to $gv_S$. \qed  

\begin{rem}{\em
%(1.) 
Iversen uses the description of spheres $S$ in $\R^n$ by the quadratic equations 
$$
b \<x, x\> - 2\<x, f\> + a=0,  
$$
Assuming the normalization $|f|^2 -ab=1$, the inversive coordinates then become (for $r<\infty$) 
$$
(b, f, a). 
$$
The natural invariant, the negative of the cosine of the angle between the spheres, is given by 
$$
\<f_1, f_2\> - \frac{1}{2}( a_1 b_2 + b_1 a_2),
$$ 
see \cite[sect. I.8]{Iversen}, i.e.  the Lorentzian inner product for vectors in $\{Q=-1\}$. In the case when $S$ is a cooriented hyperplane, $b=0$, $f$ is the unit vector normal to $S$ and $a=\hat{r}$ is the coradius. 
This gives an alternative proof of the lemma. 

%(2.) 

}
\end{rem}

Fixing one  sheet of the two-sheeted hyperboloid $\{Q=-1\}$ as a model of $\bH^{n+1}$, the original sphere $S$ corresponds 
to the boundary at infinity of the intersection with $\{Q=-1\}$ of the plane $Q$-orthogonal to $v_S$. 
Then under the isomorphism $j_*$ induced by $j$, $\mathrm{Isom}(\bH^{n+1}) \cong O^+_Q(\R)$, where $O_Q^+$ is the ``orthochronous''  
subgroup of $O_Q$ which preserves the two sheets of $Q=-1$ (rather than allowing them to interchange). 

\medskip

The pair $(\R^{n+2},Q)$ defines a hyperbolic quadratic space. Write $(\R^{n+2})^*$ for the dual vector space, and $Q^*$ for the induced dual form.
The key observation elucidating the role of isotropic vectors in the study of the arithmetic of sphere packings is the following.

\begin{lem}\label{lem:bbIs}
The  ``bend'' covector $\bb=(0,\dots,0,-2)\in(\R^{n+2})^*$ is isotropic, 
\be\label{eq:bbIso}
Q^*(\bb)=0
\ee 
and captures the bend of a sphere $S$ with inversive coordinates $v_S=(\frac1r,\frac zr, \frac1{\hat r})^t$. That is,
\be\label{eq:bbVS}
\bb(v_S)=\frac1r.
\ee
Similarly, the co-bend covector $\widehat\bb=(-2,0,\dots,0)$ is also isotropic, and has $\widehat\bb(v_S)=\frac1{\hat r}$.
In the dual inner product $\<\cdot,\cdot\>_*$ defined by $Q^*$, we have
\be\label{eq:bbbhat}
\<\bb,\widehat\bb\>_*=-2.
\ee
\end{lem}

\pf Direct and elementary computation.\epf

\medskip 
We can also identify the hyperbolic space $\H^{n+1}$ with a component of the  
two-sheeted hyperboloid $\{\cF =-1\}$, where
$(V,\cF)$ is another real hyperbolic quadratic space  of the same dimension $n+2$.
% To choose coordinates on $\bH^{n+1}$ means to make an invertible linear change of variables from $(\R^{n+2},Q)$ to another ; 
A convenient way to choose coordinates here is as follows. 
Let $\< \cdot, \cdot \>_\cF$ denote the bilinear form on $V$ corresponding to $\cF$, and let  $\cF^*$ be the dual form on the dual space $V^*$. 
%$$
%\cF^*(\al_v)= \cF(v), \hbox{~~where~~} \al_v(x)= \<v, x\>_\cF, x\in V. 
%$$
\begin{lemma}\label{lem:dualSplit}
The dual space $V^*$ admits an orthogonal splitting 
\be\label{eq:Vsplit}
V^*=V_1^*\oplus V_2^*
\ee 
so that $\dim V_1^*=2$
with
\be\label{eq:FonV1}
\bigg.\cF^*\bigg|_{V_1^*}=\mattwo 0{-2}{-2}0,
\ee
and $\cF^*$ restricts to 
a definite form on the second factor. 
\end{lemma}
\proof Take two linearly independent  light-like covectors $\alpha_0, \alpha_{n+1}\in V^*$ (i.e. $\cF^*(\alpha_0)=\cF^*(\alpha_{n+1})=0$), and rescale $\alpha_{n+1}$ to ensure that $\<\alpha_0,\alpha_{n+1}\>_{\cF^*}=-2$, as in \eqref{eq:bbbhat}. Let $V_1^*$ be the span of these vectors. 
Defining $V_2^*$  to be the orthogonal complement spanned by an arbitrary orthonormal system $\alpha_1,\dots,\alpha_n$ gives the required splitting. 
%It is exactly here that we require $\widetilde\G$ to be non-uniform, so that the form $\cF$ is isotropic, and $v_1, v_2$ exist! 
\qed

\medskip
We next transition to the number-theoretic discussion. 

\begin{add}\label{add:k'}
Suppose that $\cF$ is a hyperbolic quadratic form  on an $(n+2)$-dimensional real vector space $V$ and 
  there is a basis of $V$ with respect to which  the  form is a hyperbolic over a (totally real) field $k$. Then:
  
1. There is field extension $k'\supset k$ of degree at most two with ring of integers $\fo'$, and a choice of covectors $\alpha_0,\dots,\alpha_{n+1}$ so that the splitting \eqref{eq:Vsplit} is realized with $V_1^*$ spanned by $\alpha_0,\alpha_{n+1}$ and $V_2^*=\<\alpha_1,\dots,\alpha_{n}\>$, where $\alpha_0$ can be chosen to be defined over $\fo'$. 
%\textcolor{red}{add a statement about a unit vector $u$?} 

2. If $k=\Q$ and $\cF$ is isotropic over $\Q$, then we can take $k'=k=\Q$ and $\fo'=\Z$. Moreover, 
$V$ contains a rational unit vector $u$ (of $\cF$). 

\end{add}
\pf 1. For a general field $k$, we can find an  isotropic covector $\alpha_0$  over a suitable quadratic extension $k'$.\footnote{As an aside, recall Godement's compactness criterion, that $O_\cF(\fo)$ is non-uniform if and only if $\cF$ is isotropic over $k$, see \cite{BorelHC1962}; and furthermore, if this is the case, then $k=\Q$. Indeed, if $\cF$ represents $0$ nontrivially, then so does every Galois conjugate $\cF^\sigma$. But $\cF^\sigma$ is definite since $\cF$ is hyperbolic, so there can be no other Galois conjugates, and $k=\Q$.} We may clear denominators to ensure that $\alpha_0$ is defined over $\fo'$. Then we proceed as in the proof of \lemref{lem:dualSplit}, finding a second isotropic vector $\alpha_{n+1}$ (again over $k'$) and rescaling to ensure \eqref{eq:FonV1}. Lastly, 
choose an orthonormal basis $\al_1,\dots,\al_n$ for the orthogonal complement $V_2^*$ of $V_1^*$.

2. Since $\cF$ is isotropic over $\Q$, so is the dual form $\cF^*$. We then let $\alpha_0$, $\alpha_{n+1}$ denote linearly independent 
rational isotropic vectors of $\cF^*$. As before, by rescaling we may assume that $\alpha_0$  is an integer vector and  \eqref{eq:FonV1} is satisfied. In order to find 
a rational unit vector, we similarly find two rational isotropic vectors $v, w\in V$ such that $\<v,w\>_\cF= 1/2$. 
Then the vector $u:= v+w$ has $\<u,u\>_\cF=1$, as desired.   \epf

\subsection{Constructing (Super)Integral Kleinian Bugs/Packings}\

\pf[Proof of \thmref{thm:converseArith}]
We first prove part $(ii)$. Let $\Ga$ be a non-uniform $\Q$-arithmetic group of simplest type, commensurable to $O_\cF(\Z)$, where 
 $\cF$ is a hyperbolic quadratic form defined over $\Q$ (with respect to some basis for $V$). 
 Let $R_S\in \Ga$ be a reflection with respect to a rational  vector $u=w_S$:
  $$
 R_S:x\mapsto I-2{\<x,w_S\>\over\<w_S,w_S\>}w_S
  $$ 
  Hence, $a=\<w_S,w_S\>$ is rational. Dividing $\cF$ by $a$ we obtain a new rational hyperbolic form  $a^{-1}\cF$, with respect to which $w_S$ is a unit vector. 
  Then we  choose coordinates $\al_0,\dots,\al_{n+1}$ for $V^*$ as in \addref{add:k'}, adapted to the form $a^{-1}\cF$. 
  
Since $\Ga$ is commensurable to  $O_\cF(\Z)$, all the vectors in the $\Ga$-orbit of the vector $u$ have uniformly bounded denominators, hence, their bends  $\al_0(\G u)$ also have uniformly bounded denominators. Under the correspondence, $\gamma(S)\mapsto \gamma(u)$ defined in  \secref{sec:inversive} (with respect to the form $a^{-1}\cF$), 
the radius of the sphere $\gamma(S)$ equals $\alpha_0(u)$. Hence, by a suitable rescaling in the Euclidean space $\R^n$ (equivalently, rescaling the choice of $\al_0$ to clear the denominators), all the bends can be made integral, as claimed. 

The proof of part $(i)$ is now clear. By rescaling the form $\cF$, we get a new form  rational hyperbolic form  $a^{-1}\cF$ such that the reflection $R_S\in \G$ is given by a 
unit (with respect to $a^{-1}\cF$) vector $w_S$. The group $\widetilde\Ga$, of course, preserves the new form. 
Since, by the assumption, the action of the symmetry group $\G_S< \widetilde\Ga$ on the packing $\sP$ is transitive, 
so is the action of  $\widetilde\Ga$ on the superpacking $\widetilde\sP$. Thus, by Part (ii), 
after a suitable rescaling, the superpacking $\widetilde\sP$ becomes integral. 
\epf

\bigskip 
We now prepare for the proof of  \thmref{thm:quadExt}. Before getting to questions of (super)integrality, we construct the necessary packings as follows.

\begin{prop}\label{prop:Millson}
Let  $O_{{\mathcal F}}(\fo)<\mathrm{Isom}(\bH^{n+1})$ 
be a $k$-arithmetic lattice of simplest type. 
Then there is a sequence of conformally inequivalent Kleinian packings $\sP_j$ such that 
$O_{{\mathcal F}}(\fo)$ is commensurable to a supergroup $\widetilde{\Ga}_j$ of 
 $\sP_j$ which acts transitively on $\sP_j$. 
\end{prop}
\proof 
%In fact, we will construct an infinite sequence $\sP_j$ of Kleinian packings whose supergroups are commensurable to $O_{{\mathcal F}}(\fo)$. 
While this result is stated in arithmetic terms, most of the proof is non-arithmetic. 

\medskip 
{\bf Geometric setup.} Suppose that $\Ga< G=\mathrm{Isom}(\H^{n+1})$ 
is a torsion-free lattice such that the hyperbolic manifold $M=\H^{n+1}/\Ga$  
contains a properly embedded complete connected nonseparating totally 
geodesic hypersurface of finite $n$-dimensional volume $N$. The hypersurface 
$N$ lifts to a hyperplane $H$ in $\H^{n+1}$. We will assume that the reflection 
$R_H$ in $H$ normalizes the group $\Ga$ and let 
$\widetilde\Ga$ denote the subgroup of $G$ generated by $R_H$ and $\Ga$. 
The group $\widetilde\Ga$ contains $\Ga$ as an index 2 subgroup; hence, 
$\widetilde\Ga$ is again a lattice. Since $R_H$ normalizes $\Ga$, it descends to an isometric 
reflection $\si: M\to M$ fixing $N$ pointwise. 

\begin{rem}{\em
Abundance of examples of this type (comping from arithmetic groups) 
was first established by John Millson in \cite{Millson}. More precisely, he proved that every 
$k$-arithmetic lattice $O_{{\mathcal F}}(\fo)$ in $O({n+1,1})$ 
of simplest type (where $\fo$ is the ring of integers of the field $k$ and ${\mathcal F}$ 
is a hyperbolic quadratic form over $k$) is commensurable to a group $\widetilde\Ga$ as above. }
\end{rem}

Our goal is define a sequence of Kleinian packings $\sP_j$ with the supergroups $\widetilde\Ga_j$ 
commensurable to $\widetilde\Ga$ and acting transitively on $\sP_j$.  

%which is, moreover, superintegral if $k=\Q$ and  ${\mathcal F}$ is diagonal. In fact, we will also construct a sequence of packings whose supergroups $\widetilde\Ga_k$ are finite-index subgroups of $\widetilde\Ga$. (This sequence will be used in section \ref{}.) 

The hypersurface $N$ does not separate $M$ and, hence, defines a nontrivial element $\xi$ of $H^1(M)$, which is  
Poincar\'e-dual to the locally-finite fundamental class of the hypersurface $N$. Since $\si$ fixes $N$ pointwise, $\si^*(\xi)=-\xi$.  
The class $\xi$ defines a homomorphism
$$
\phi: \Ga\to H_1(M)\to \Z. 
$$
Let $p: \hat{M}\to M$ denote the infinite cyclic covering corresponding to the kernel $Ker(\phi)$ of $\phi$. 
 Since  $\si^*(\xi)=-\xi$, $Ker(\phi)$ is $R_H$-invariant; hence $\si$ lifts to a reflection $\tau: \hat{M}\to \hat{M}$ fixing 
 pointwise one of the components $N_0$ of the preimage of $N$ in $\hat{M}$. We let $D$ denote a component 
 of $p^{-1}(M - N)$ whose boundary contains $N_0$. Then $\bar{D}$ is a fundamental domain of the  action of 
 the deck-transformation group $\Z=\<\theta\>$  of the regular covering $p: \hat{M}\to M$. 
Furthermore,  for each $j\in \N$,  
$$
D_j:= \bigcup_{-j\le i\le j} \theta^i(\bar{D})=  \bigcup_{0\le i\le j} \theta^i(\bar{D}) \cup \tau(  \bigcup_{0\le i\le j} \theta^i(\bar{D})) 
$$
is a fundamental domain for the index $2j$ subgroup in $\Z$. Each domain $D_j$ is $\tau$-invariant, 
has finite volume ($2j$ times the volume of $M$) and two boundary component, 
both totally geodesic in $\hat{M}$ and isometric to $N$ via the restriction of the covering map $p$.  

\begin{rem}\label{rem:exhaust} 
{\em
For future reference, we record the following obvious properties of the domains $D_j$: 

1. $cl(D_j)\subset D_{j+1}$ for each $j$. 

2. $\bigcup_{j\ge 1} D_j= \hat{M}$. 
}
\end{rem}

\begin{figure}%[htbp] %  figure placement: here, top, bottom, or page
   \centering
\includegraphics[width=.7\textwidth]{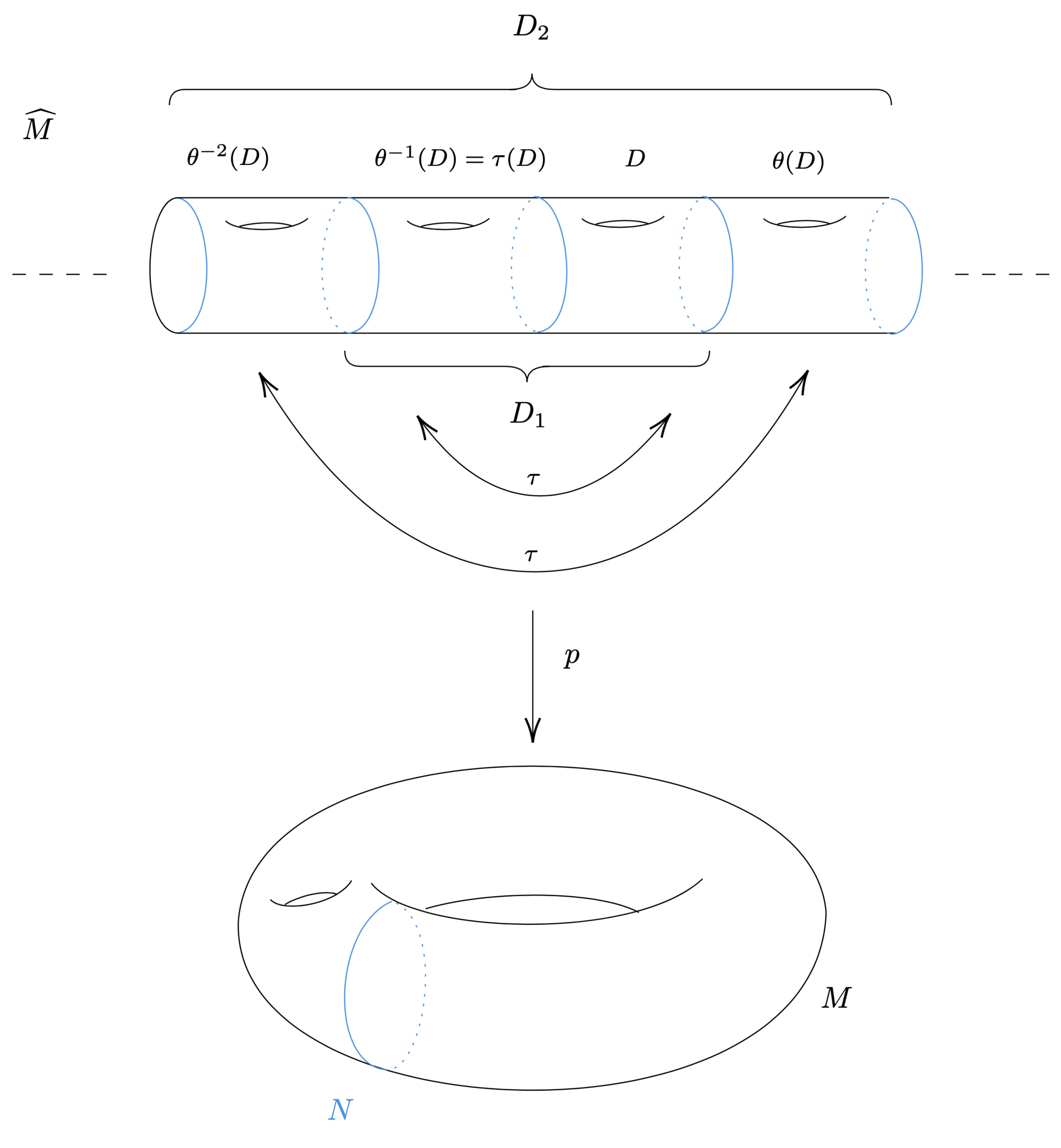}
\caption[Caption for LOF]{Infinite cyclic covering\protect\footnotemark}
%             \caption{Infinite cyclic covering.\footnote{.} }
   \label{fig:covering}
\end{figure}

The fundamental group $\Ga_j$ of $D_j$ embeds in $\pi_1(\hat{M})$ (since $D_j$ has totally-geodesic boundary) and, hence, $\pi_1(M)$. Since $\pi_1(M)$ is isomorphic to the lattice $\Ga$, the group $\Ga_j$ will be identified with a subgroup (again denoted $\Ga_j$) 
of $\Ga$. Then the preimage of $D_j$ in $\H^{n+1}$ under the covering map $q: \H^{n+1}\to \hat{M}$ contains a (unique) 
$\Ga_j$-invariant component $P_j=\tilde{D}_j$. Since $\partial N$ is totally-geodesic, $P_j$ is a convex domain  with totally-geodesic boundary. By the construction, each boundary component of $P_j$ is fixed by a reflection in $\widetilde\Ga$ conjugate to $R_H$. Thus, the collection of hyperplanes bounding $P_j$ defines a sphere packing $\sP_j$ with the reflection group 
$\Ga_{R,j}$, generated by reflections in the walls bounding $P_j$. The group $\Ga_j$ acts on $P_j$ with 
quotient of finite volume (isometric to $D_j$); in particular, the limit set $\La_j$ 
of $\Ga_j$ equals $\geo P_j$ and $P_j$ is the closed convex hull of this limit set, $P_j= C_{\Ga_j}$. We conclude that 
$\Ga_j$ is geometrically finite. Thus, $\sP_j$ is a Kleinian packing.  

\footnotetext{Many thanks to Jules Flin for this drawing.}

The supergroup of the packing generated by $\Ga_{R,j}$ and $\Ga_j$, however, acts on 
$\sP_j$ with {\em two orbits}: These two orbits correspond to the two connected components of $\partial D_j$.  
To fix this problem, we note that the reflection $\tau$ discussed above swaps these boundary components. 
Hence, we lift $\tau$ to a reflection $\tilde\tau$ in $\H^{n+1}$ preserving $P_j$ and let $\widetilde\Ga_j$ denote the subgroup 
of $\Isom(\H^{n+1})$ generated by $\Ga_j$ and $\tilde\tau$. The group $\widetilde\Ga_j$, which is an index 2 extension of 
$\Ga_j$, is still a symmetry group of the packing $\sP_j$ and acts transitively on the packing. 

We are almost done with the proof of the proposition. It remains to show that the packings $\sP_j$ are conformally inequivalent for different $j$'s; equivalently, we claim that the convex subsets $P_j\subset \H^{n+1}$ are pairwise non-isometric. 

 To distinguish the sets $P_j$, we define the following invariant:
$$
\rho_j= \sup_{\bar{x}\in D_j} d(\bar{x}, \partial D_j)= 
\sup_{{x}\in P_j} d({x}, \partial P_j). 
$$ 
Here $d({x}, \partial P_j)$ denotes the minimal distance from $x$ to the points of $\partial P_j$, similar for 
$d(\bar{x}, \partial D_j)$. Since each $P_j$ has finite volume, $\rho_j$ is finite for each $j$. 
Then  \rmkref{rem:exhaust}(1) implies that the sequence $\rho_j$ is strictly  increasing with $j$ and, hence $P_j$'s are pairwise non-isometric.  \qed

\begin{rem}
{\em
1. In  \secref{sec:Hdim} we will give a different argument for conformal inequivalence of packings (possibly after passing to a subsequence) using Hausdorff dimensions of the limit sets of the groups $\Ga_j$. 

2. Two Kleinian packings are conformally equivalent if and only if they are quasiconformally (more precisely, quasisymmetrically) equivalent, provided that $n\ge 2$, see \cite{Frigerio2, BKM}. Hence, we obtain infinitely many quasiconformally  inequivalent packings in every dimension $n\ge 2$. 
}
\end{rem}

\pf[Proof of \thmref{thm:quadExt}]\

We are given $O_\cF(\fo)$, a $k$-arithmetic hyperbolic lattice of simplest type, with $\fo$ the ring of integers of $k$. 
In \propref{prop:Millson}, we constructed a Kleinian packing $\sP=\G_S\cdot S_0$ with a symmetry group $\Ga_S$ 
which acts transitively on the spheres in the packing, and an arithmetic super-symmetry group $\widetilde\Ga$ 
(commensurable to $O_\cF(\fo)$). The group $\widetilde\Ga$  contains a reflection $R$ through the sphere $S_0$. As in the proof of \thmref{thm:converseArith}, the normal vector $w\in V$ to $S_0$ has all coordinates in $k$. By \addref{add:k'}, 
after a conformal change of coordinates on $\R^n$,
the ``bend'' covector $\al_0$ is defined over $\fo'$, the ring of integers of a quadratic extension $k'$ of $k$. The supergroup $\widetilde\G$ of the packing  is commensurable to $O_\cF(\fo)$,  and so the orbit of $w$ under the supergroup is defined entirely over $k$; it is only when we measure the bends using the covector $\al_0$ (after a suitable rescaling) that we obtain elements of $\fo'$. Regardless, the superpacking has all bends in $\fo'$, as claimed. \epf

\pf[Proof of the  Classification \thmref{thm:class}]\

The same proof gives the forward direction of the Classification \thmref{thm:class}; indeed, if $\cF$ is defined over $\Q$ and isotropic, then by the same argument as above, the packing constructed in \propref{prop:Millson} is superintegral by Part $(i)$ of \thmref{thm:converseArith}. 
The backwards direction is a direct consequence of the Subarithmeticity \thmref{thm:arithmBugs}, which we turn to now. 
We will give two proofs, one that is basically identical to the proof of \cite[Thm 19]{KN}, and another using somewhat different ideas. 
\epf

\pf[First Proof of \thmref{thm:arithmBugs}]\

Let $Q$ be the standard quadratic form of the signature $(n+1,1)$ as in equation \eqref{eq:Qdef},
 $\G<G=O_Q^+(\R)$
be a discrete, Zariski dense subgroup
 acting on the inversive coordinates $v_{S_0}$ of a sphere $S_0$, so that
 the bends, that is, first entries, in the orbit $\cO=\G v_{S_0}$ are all integers. 
 The action of $\G$ is on the left (on column vectors) and involves all the entries of $v_{S_0}$; we conjugate it to a right action just on (row vectors of) bends, as follows.
 
 By the Zariski density of $\G$, the orbit $\cO$
 is also Zariski dense in the one-sheeted hyperboloid $Q=1$, and hence 
contains $n+2$ linearly independent vectors $\{v_1=v_{S_0},v_2,\dots,v_{n+2}\}\subset\cO$. 
These vectors provide a coordinate system for the dual vector space $(\R^{n+2})^*$ so that if one applies these coordinates to the covector $\al_0$ (that is, the first coordinate in $\R^{n+2}$), one gets integers.
Moreover, for each $\gamma\in \Gamma$, the pairing of vectors and covectors $\<\gamma v_i, \al_0\>=\<v_i, {}^*\gamma^{-1}\al_0\>$ is also an integer. Hence, the vectors $v_i$ evaluate to integers on the $\Gamma$-orbit of the covector $\al_0$ in   $(\R^{n+2})^*$. Hence $\Gamma$ preserves the finite index sublattice $L$ generated by the $\Gamma$-orbit of $\al_0$ in  $(\Z^{n+2})^*$; here  $(\Z^{n+2})^*$ consists of covectors with integer coordinates with respect to the coordinate system given by $v_1,..., v_{n+2}$.
Thus $\Gamma<G^L$ is a subgroup of $G^L=\{g\in O_Q^+(\R):gL=L\}$.
The latter group is easily seen to be (and is sometimes taken to be the very definition of) a congruence subgroup.
Note also that the co-vector $\al_0$ is rational with respect to the integral structure given by $L$ and isotropic,  see  \lemref{lem:bbIs}.  It follows that $G^L$, regarded as a lattice acting on the 
dual vector space 
$V^*$, is non-uniform, i.e. contains a unipotent element. Then $G^L$ itself, regarded as a subgroup of $G$, contains a unipotent element, and, hence, is non-uniform. To conclude:    $\G$ is a subgroup of a non-uniform $\Q$-arithmetic hyperbolic group of simplest type, as claimed. 
\epf

\pf[Second Proof of \thmref{thm:arithmBugs}]\

Let
$\G<O_Q^+(\R)$ 
be a discrete, Zariski dense  subgroup
 acting on the inversive coordinates $v_{S_0}$ of a sphere $S_0$, so that
 the bends, that is, first entries, in the orbit $\cO=\G v_{S_0}$ are all integers. 
 The first coordinate on $\R^{n+2}$ is a (nonzero) linear functional $\al$ on the real vector space $V= \R^{n+2}$. The key is the following general lemma:
 
 \begin{lemma}
 Let $V$ be a finite-dimensional real vector space,  $\Gamma< GL(V)$ an irreducible subgroup, i.e.  a subgroup which has no proper invariant subspaces. 
 Let $v\in V, \al\in V^*$ be nonzero vectors with the property that $\al(g v)\in \Z$ for all $g\in \Gamma$. Define the $\Z$-submodule $L$ in $V$ generated by the orbit $\Gamma\cdot v$. 
 Then $L$ is a free $\Z$-module of rank equal to the dimension of $V$. 
 \end{lemma}
 \proof The group $\Gamma$ obviously acts on $L$ by automorphisms and $\al$ still takes only integer values on $L$. The irreducibility of $\Ga$ implies that $L$ spans $V$ as a real vector space. It remains to prove that $L$ is a discrete subgroup of $V$ regarded as an abelian Lie group. Let $\bar{L}$ denote the closure of $L$ in the classical topology on $V$ and let 
 $W:=\bar{L}_0$ be the identity component of this Lie subgroup of $V$; this component is a (real) linear subspace in $V$. The group $\Gamma$ preserves this subspace. In view of irreducibility of $\Ga$, the subspace $W$ is either $\{0\}$ or the entire $V$. However, $\al\ne 0$ still takes only integer values on $W$, hence, $W\ne V$ and we conclude that $W=\{0\}$,  i.e. $L$ is a discrete subgroup of $V$.  
 \qed 

\medskip 
 We apply this lemma in our setting. The submodule $L\subset V=\R^{n+2}$ defines an integral structure on $V$. We claim that the quadratic form $Q$ is rational with respect to this integral structure. Since $G_{\mathbb Z}$, the set of integer points in $G=O_Q$, contains $\Ga$, it is Zariski dense in $G$, hence, is an (arithmetic) lattice. It follows from Exercise 
 4 in \cite[Sect. 5A]{Witte} that $G$ is defined over $\Q$. We claim that the form $Q$ is also defined over $\Q$. The proof is the same as the one of 
 Exercise  4 in \cite[Sect. 5A]{Witte}: Consider the  vector space $U$ of all quadratic forms on $V$. In view of Zariski density of $\Ga$ in $G$, there is a unique $\Ga$-invariant line 
 in $U\otimes \C$. Since $\Ga$ consists of integer matrices (with respect to the integer structure on $V$ defined by $L$), for every Galois automorphism $\sigma$ of $\C$, for 
 every $q\in U$, we have
 $$
 \ga^*(q)^{\sigma}= \ga^*(q^\sigma),\ga\in \Ga. 
 $$
Therefore, since $Q$ is $\Ga$-invariant, so are the forms $Q^\sigma, \sigma\in Gal(\C)$. Thus,  the forms $Q^\sigma$ belong to the line $\C Q$, i.e. for every $\sigma\in Gal(\Q)$ there exists $z\in \C^\times$ such that 
$$
Q^\sigma =z Q. 
$$ 
 Lastly, $Q(v)=1$ and $v\in L$, hence, 
 $$
z= z Q(v)= Q^\si(v)= Q^\si(v^\si)= (Q(v))^\si=1. 
$$
It follows that   $Q$ itself is a rational form with respect to the rational structure on $V$ defined by the lattice $L$. 

Lastly, the argument that $\Ga$ is contained in a nonuniform lattice is the same as in the first proof: The 
covector $\alpha_0$ is isotropic and rational with respect to the integral structure defined by the lattice $L$.  
\epf

\subsection{Caveats and Examples}\label{sec:exs}\

We collect here some examples that illustrate various caveats given in the introduction to the main theorems.
We begin with an explicit example of a superintegral Kleinian packing which is not crystallographic.

\begin{figure}
\includegraphics[width=2in]{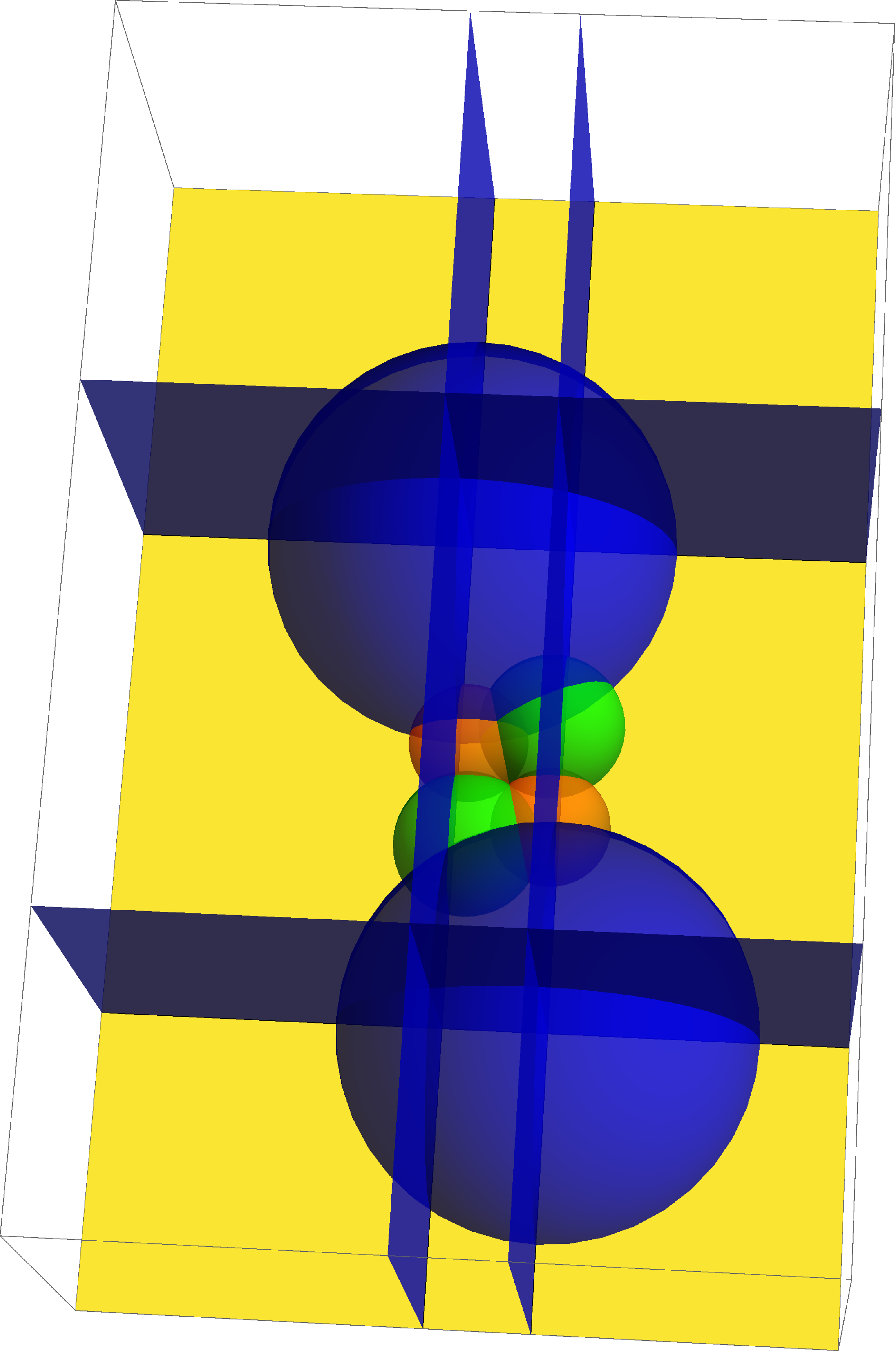}
\
\includegraphics[width=2.5in]{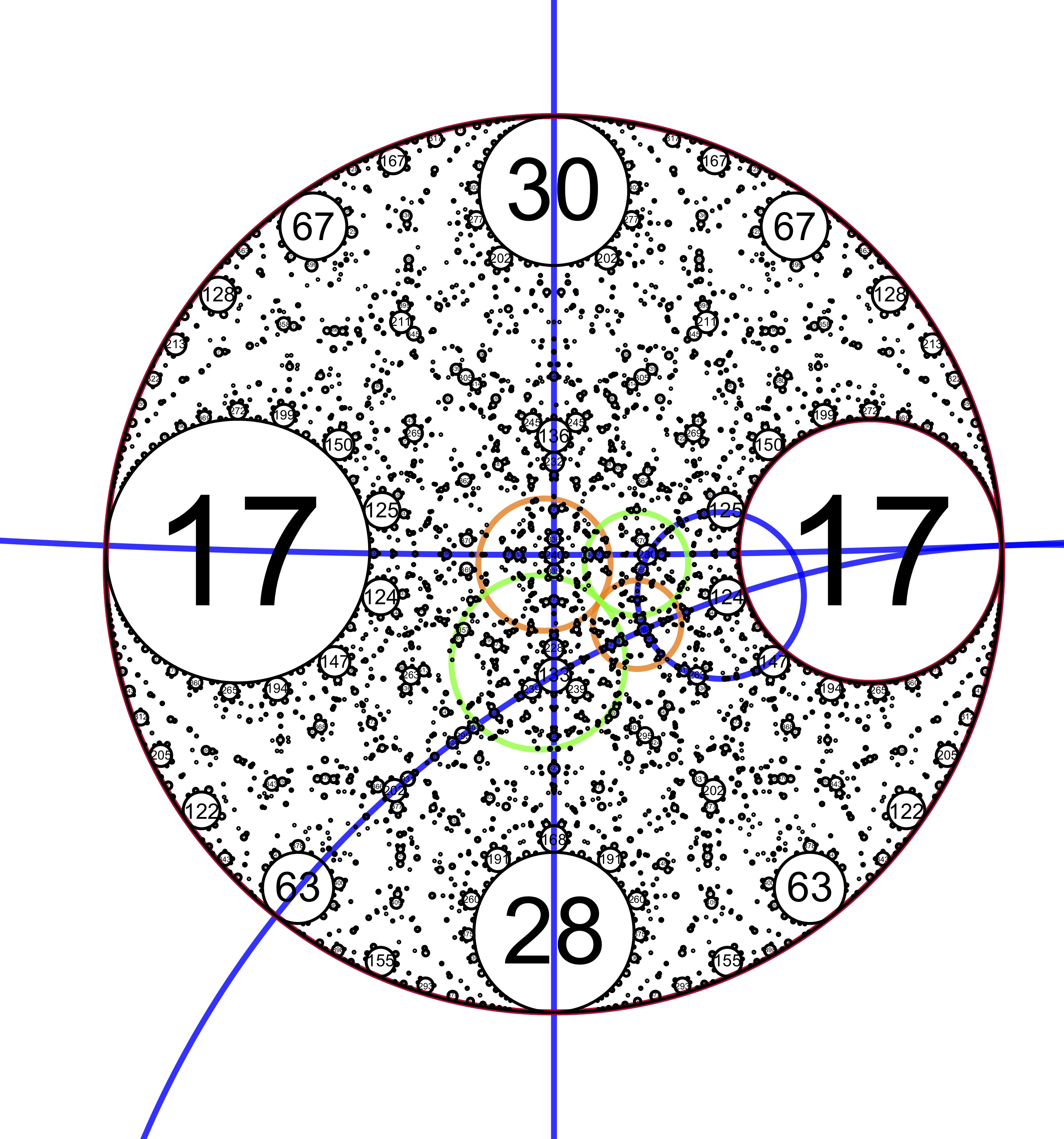}

$(a)$ \hskip 1in $(b)$
\caption{$(a)$ A fundamental domain for the extended Bianchi group $\widehat{Bi}(23)$, and $(b)$ a  superintegral Kleinian packing attached to it.}
\label{fig:Bianchi}
\end{figure}

\begin{example}\label{lem:notCryst}
{\em
The extended Bianchi group $\widehat{Bi}(23)$ (see footnote \ref{foot:Bianchi}) %containing the group $\PSL_2(\cO_{23})$ 
is not reflective (see \cite{BelolipetskyMcleod2013}). %; here $\cO_{23}=\Z[{1+\sqrt{-23}\over2}]$ is the ring of integers of $\Q(\sqrt{-23})$. 
Indeed, applying Vinberg's algorithm \cite{Vinberg1972} to the quadratic form $f=-2xy+2z^2+2zw+\frac{23+1}2w^2
%-2x_1x_2+2x_3^2+2x_3x_4+\frac{23+1}2x_4^2
$ shows that the subgroup of $O_f(\Z)$ generated by all reflections has infinite index in $O_f(\Z)$. One can give a  fundamental domain in $\bH^3$ for $O_f(\Z)$ as shown in \figref{fig:Bianchi}$(a)$; here the blue walls act by reflections,  %this is a bit strange terminology, but its clear what it means
and there is a pair of commuting unipotent elements, one  identifying the green walls, and another identifying the orange ones. By Part 2(ii) of the Structure \thmref{thm:lattice}, we can construct a superintegral packing from this fundamental polyhedron; see \figref{fig:Bianchi}$(b)$.
So this is a superintegral packing which is Kleinian but not crystallographic.\footnote{Note 
that this particular packing appeared previously in work of Stange  \cite{Stange2018} and Martin \cite{Martin2019} but was not recognized there as being dense (and hence wasn't considered a ``packing'' by our definition), due to the significant distances between disjoint circles. The general theory given here makes this density apparent.  %transparent.
}
}
\end{example}

As mentioned above \rmkref{rmk:10}, a superintegral Kleinian packing can have all spheres disjoint (as observed already in \cite{KN}); but its supergroup must still be non-uniform and $\Q$-arithmetic. One such is the following.

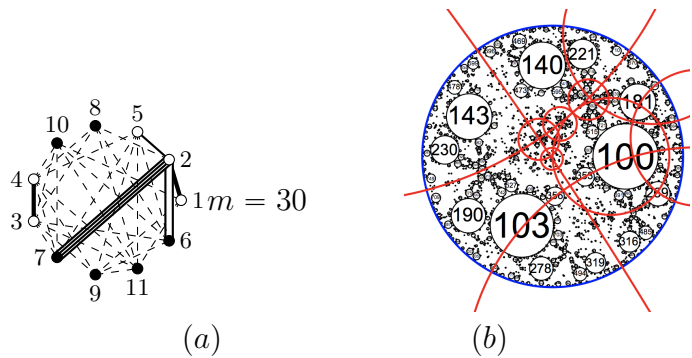
\begin{figure}
 \begin{tikzpicture}
    \coordinate (one) at (1, 0);
    \coordinate (two) at (0.84,0.54);
    \coordinate (five) at (0.42,0.91);
    \coordinate (eight) at (-0.14,0.99);
    \coordinate (ten) at (-0.65,0.76);
    \coordinate (four) at (-0.96,0.28);
    \coordinate (three) at (-0.96,-0.28);
    \coordinate (seven) at (-0.65,-0.76);
    \coordinate (nine) at (-0.14,-0.99);
    \coordinate (eleven) at (0.42,-0.91);
    \coordinate (six) at (0.84,-0.54);

    \draw[ultra thick] (one) -- (two);
    \draw[dashed] (one) -- (nine);
    \draw[dashed] (one) -- (ten);
    \draw[dashed] (one) -- (eleven);
    \draw[thick] (two) -- (five);
    \draw[thick, double distance = 2pt] (two) -- (six);
    \draw[thick, double distance = 2.5pt] (two) -- (seven); 6.0
    \draw[thick, double distance = 0.3pt] (two) -- (seven); 6.0
    \draw[dashed] (two) -- (eight);
    \draw[ultra thick] (three) -- (four);
    \draw[dashed] (three) -- (six);
    \draw[dashed] (three) -- (seven);
    \draw[dashed] (three) -- (eight);
    \draw[dashed] (three) -- (nine);
    \draw[dashed] (three) -- (ten);
    \draw[dashed] (three) -- (eleven);
    \draw[dashed] (four) -- (five);
    \draw[dashed] (four) -- (seven);
    \draw[dashed] (four) -- (eight);
    \draw[dashed] (four) -- (nine);
    \draw[dashed] (four) -- (ten);
    \draw[dashed] (four) -- (eleven);
    \draw[dashed] (five) -- (six);
    \draw[dashed] (five) -- (seven);
    \draw[dashed] (five) -- (nine);
    \draw[dashed] (five) -- (eleven);
    \draw[dashed] (six) -- (eight);
    \draw[dashed] (six) -- (nine);
    \draw[dashed] (six) -- (ten);
    \draw[dashed] (seven) -- (ten);
    \draw[dashed] (eight) -- (eleven);
    \draw[dashed] (nine) -- (eleven);
    \draw[dashed] (ten) -- (eleven);

    \filldraw[fill=white] (one) circle (2pt) node[right] {\scriptsize 1} node[right, xshift = 0.2cm] {$m=30$};
    \filldraw[fill=white] (two) circle (2pt) node[right] {\scriptsize 2};
    \filldraw[fill=white] (three) circle (2pt) node[left] {\scriptsize 3} ;
    \filldraw[fill=white] (four) circle (2pt) node[left] {\scriptsize 4};
    \filldraw[fill=white] (five) circle (2pt) node[above] {\scriptsize 5};
    \filldraw (six) circle (2pt) node[right] {\scriptsize 6};
    \filldraw (seven) circle (2pt) node[left] {\scriptsize 7};
    \filldraw (eight) circle (2pt) node[above] {\scriptsize 8};
    \filldraw (nine) circle (2pt) node[below] {\scriptsize 9};
    \filldraw (ten) circle (2pt) node[above] {\scriptsize 10};
    \filldraw (eleven) circle (2pt) node[below] {\scriptsize 11};
  \end{tikzpicture}
\hskip.4in
    \includegraphics[width=.3\textwidth]{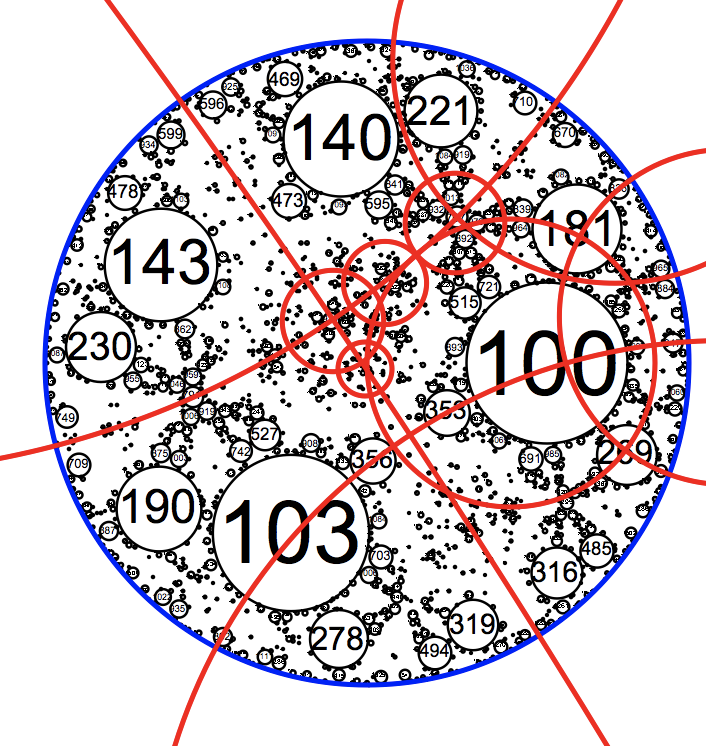}
    
$(a)$\hskip1.3in$(b)$  
\caption{$(a)$ The Coxeter diagram for $\widehat{Bi}(30)$, and $(b)$ a superintegral packing with all spheres disjoint.}
\label{fig:nonT}
\end{figure}

\begin{example}\label{ex:nonTangent}
{\em
Consider the extended Bianchi group $\widetilde\G=\widehat{Bi}(30)$.
%to be the 
%maximal discrete group of $\mathrm{Isom}(\bH^3)$ 
%containing $\PSL_2(\cO_{30})$, with $\cO_{30}$ the ring of integers of $\Q(\sqrt{-30})$.
Equivalently, $\widetilde\G\cong O_\cF(\Z)$ is %isomorphic to 
the integer orthogonal group preserving the form $\cF=-xy+z^2+30w^2$.
This group is reflective (see  \cite{BelolipetskyMcleod2013}), and 
applying Vinberg's algorithm \cite{Vinberg1972}
gives the Coxeter diagram%
\footnote{%
The diagram has nodes for each facet, and facets that are orthogonal are not connected; the dihedral angle $\pi/3$ is denoted by a single line, angle $\pi/4$ is a double line, and angle $\pi/\infty=0$ is a thick solid line. Nodes of separated facets are connected with a dotted line.
}  
shown in \figref{fig:nonT}. The node marked ``8'' is totally isolated from the others, being either orthogonal or some given distance apart from the other generating reflective walls. (So are nodes ``9'', ``10'', and ``11'' but we choose to use ``8''.) Dropping this wall from the generators and letting the remaining generators act on it by reflections (as in the Structure \thmref{thm:lattice}) gives the superintegral crystallographic packing shown in \figref{fig:nonT}(b).
}
\end{example}

\medskip
As observed already in \cite{KN}, if 
%in the above construction we replace the normal closure of a single reflection $R\in \widetilde{\Ga}$ by the normal closure of two non-conjugate reflections, 
a packing cannot be realized as the orbit of a single sphere (that is, there is no symmetry group for which the action on the spheres in the packing is transitive),
then the resulting  packing (or bug) need {\it not} be superintegral (or even integral); cf. \rmkref{rmk:nonInt}. Here is an explicit example of a non-integral bug with non-uniform, $\Q$-arithmetic supergroup.

\begin{figure}
\includegraphics[width=.3\textwidth]{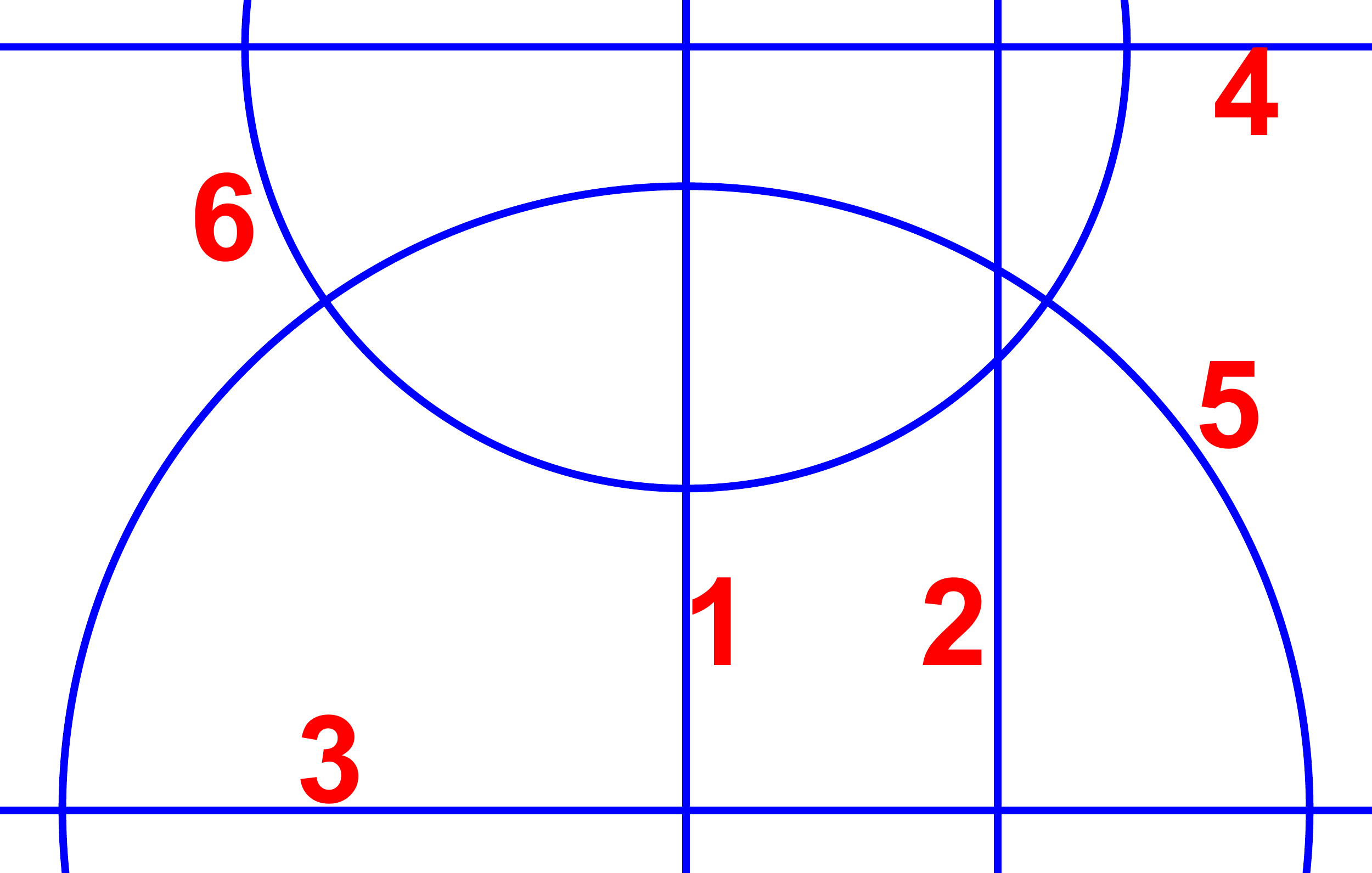}\hskip.4in
 \begin{tikzpicture}
    \coordinate (one) at (-2, 0);
    \coordinate (two) at (-1, 0);
    \coordinate (three) at (-1, 1);
    \coordinate (four) at (0, 1);
    \coordinate (five) at (0, 0);
    \coordinate (six) at (-2, 1);

    \draw[ultra thick] (one) -- (two);
    \draw[thick] (two) -- (five);
    \draw[thick, double distance = 2pt] (two) -- (six);
    \draw[ultra thick] (three) -- (four);
    \draw[dashed] (three) -- (six);
    \draw[dashed] (four) -- (five);

    \filldraw[fill=white] (one) circle (2pt) node[below] {\scriptsize 1};
    \filldraw[fill=white] (two) circle (2pt) node[below] {\scriptsize 2};
    \filldraw[fill=white] (three) circle (2pt) node[above] {\scriptsize 3};
    \filldraw[fill=white] (four) circle (2pt) node[above] {\scriptsize 4};
    \filldraw[fill=white] (five) circle (2pt) node[below] {\scriptsize 5};
    \filldraw (six) circle (2pt) node[above] {\scriptsize 6};
  \end{tikzpicture}
  
$(a)$\hskip1.3in$(b)$  
\caption{$(a)$ The spheres with inversive coordinates \eqref{eq:VinSph} in  \exref{ex:nonInt}. $(b)$ Their Coxeter diagram.}
\label{fig:Cox}
\end{figure}

\begin{example}\label{ex:nonInt}
\emph{
Let 
 $\widetilde\G$
 be the
 extended Bianchi group 
  $\widetilde\G=\widehat{Bi}(6)$.
   %is the maximal discrete group of $\mathrm{Isom}(\bH^3)$ containing $\PSL_2(\cO_6)$, with $\cO_6$ the ring of integers of $\Q(\sqrt{-6})$. 
   Equivalently, $\widetilde\G\cong O_\cF(\Z)$ is isomorphic to the integer orthogonal group preserving the form $\cF=-xy+z^2+6w^2$.
That it is reflective (generated by reflections) is essentially due to Bianchi  \cite{Bianchi1892}; see also  \cite{BelolipetskyMcleod2013}. % for a complete classification of reflective extended Bianchi groups. 
Applying Vinberg's algorithm % from \cite{Vinberg1972}
% \textcolor{red}{this needs a further explanation: I think the point is that this algorithm is searching for generators of the reflection group, more precisely, for the space-like vectors such that reflections in their orthogonal comlements generate the reflection group; also, it would be good to give a reference to he algorithm} 
%(see, e.g., the exposition in \cite{Bel}, or for this particular case, \cite{McLeodThesis})
to this group produces the following normal vectors
$$
({0, 0,-1, 0})^t,
({1, 0, 1, 0})^t,
({0, 0, 0, -1})^t,$$
$$ 
({6, 0, 0, 1})^t,
({-1, 1, 0, 0})^t,
({2, 2, 0, 1})^t.
$$
Making a (choice of) change of variables from $\cF$ to the universal form $Q$ in \eqref{eq:Qdef}, these correspond to spheres with the following (realization of) inversive coordinates:
$$
v_1=(0 , -1 , 0, 0 )^t, \ 
v_2=(0 , 1 , 0, 1)^t, \ 
v_3=(0 , 0 , -1, 0)^t, \ 
$$
\be\label{eq:VinSph}
v_4=(0 , 0 , 1, \sqrt{6})^t, \ 
v_5=(1 , 0 , 0, -1)^t, \ 
v_6=(\sqrt{2} , 0 , \sqrt{3}, \sqrt{2})^t.
\ee
See \figref{fig:Cox}(a) for the corresponding spheres (as circles in $\R^2$).
Writing $V$ for the $4\times 6$ matrix whose columns are $v_j$, we can compute the $6\times6$ Gramian of all $Q$-inner products:
$$
V^t Q V\ =\ \cG =
\left(
\begin{array}{cccccc}
 1 & -1 & 0 & 0 & 0 & 0 \\
 -1 & 1 & 0 & 0 & -\frac{1}{2} & -\frac{1}{\sqrt{2}} \\
 0 & 0 & 1 & -1 & 0 & -\sqrt{3} \\
 0 & 0 & -1 & 1 & -\sqrt{\frac{3}{2}} & 0 \\
 0 & -\frac{1}{2} & 0 & -\sqrt{\frac{3}{2}} & 1 & 0 \\
 0 & -\frac{1}{\sqrt{2}} & -\sqrt{3} & 0 & 0 & 1 \\
\end{array}
\right)
.
$$
%That is, when spheres corresponding to $v_i$, $v_j$ intersect with dihedral angle $\pi/k$, their $Q$-inner product is $<v_i,v_j>=-\cos\frac\pi k$. 
Equivalently, $\widetilde\G$ has the Coxeter diagram
 given in \figref{fig:Cox}(b). 
A different realization of $v_j$ from \eqref{eq:VinSph} will of course have different inversive coordinates $V$ but the  Gramian and Coxeter diagram are invariants.
}
\end{example}

\begin{figure}
\includegraphics[width=.5\textwidth]{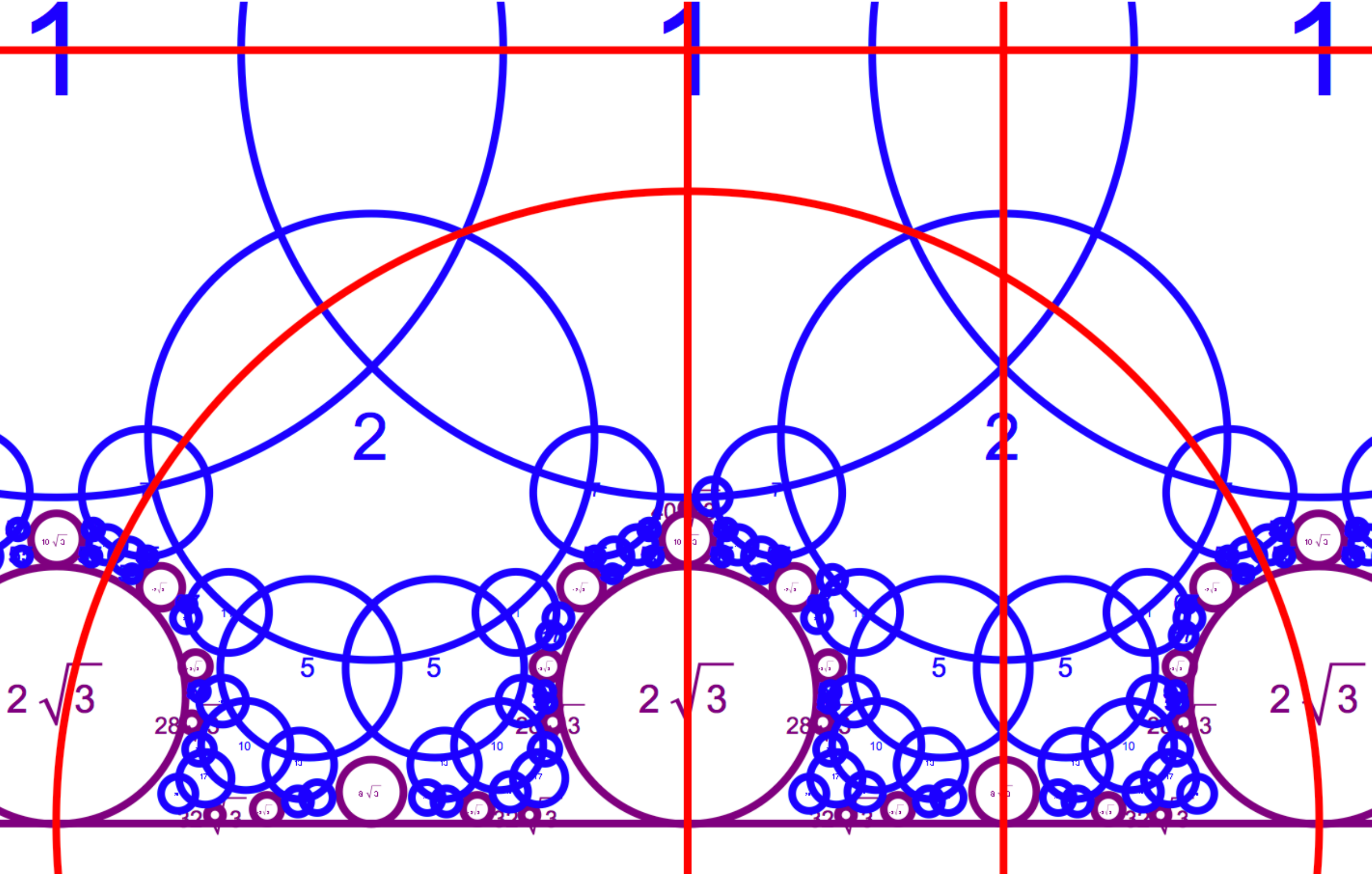}
\caption{A non-integral bug; the bend of a circle is shown at its center.}
\label{fig:nonInt}
\end{figure}

Now we construct the bug. 
%By abuse of notation, we w
Write $R_j$ for the 
 M\"obius action of reflection through sphere $v_j$, that is,
 $$
 R_j=I_{4\times4}+2 v_j\cdot v_j^t\cdot Q.
 $$
We apply Part 2(i) of the Structure \thmref{thm:lattice} with $S'=\{R_1,\dots,R_6\}$, $\<S'\>=\widetilde\G$,  and take $R=\{R_3,R_6\}$ so that $S=S'\setminus R$. The bug we obtain is then the orbit
\be\label{eq:nonIntBug}
\sB=\<S\>\cdot \{v_3,v_6\},
\ee 
as shown in 
\figref{fig:nonInt}. This particular realization of the bug is evidently non-integral, but we have not yet  ruled out that there cannot be some other conformally equivalent realization of this bug which is integral.  (Indeed, there exist realizations of the classical Apollonian packing that are non-integral.)

\begin{lemma}
There does not exist a conformal realization of the bug in \eqref{eq:nonIntBug} that is integral.
\end{lemma}
\pf
We mimic the first proof of \thmref{thm:arithmBugs} by looking at a right ``bends'' action. This time, we will do it using an over-determined system of equations, leading to an irrational linear relation among the bends, which is realization-independent (and hence the bug cannot ever be integral).

To begin, notice that the spheres 
$$v_3, ~v_6, ~~R_2\cdot v_6, ~~R_1\circ R_2\cdot v_6,  ~~\hbox{and} ~~ R_5\circ R_1\circ R_2  \cdot v_6$$ 
 are all in the bug $\sB$. Write $W$ for the $4\times 5$ matrix whose columns are 
the inversive coordinates of these spheres.
The kernel of $W$,
$$
\cK \ := \ \ker(W) \ := \ \{g\in Mat_{5\times 5}(\R):Wg=0\},
$$
is independent of the realization of the bug. Indeed,
 if $g\in\mathrm{Isom}(\bH^3)$ is any  isometry and we move the whole bug by left-acting by $g$, then the inversive coordinates matrix $W$ changes to $gW$, leaving the kernel $\cK$ invariant.
 Notice that the kernel contains, e.g.,
 $$
 K\ = \ 
\left(
\begin{array}{ccccc}
 -\sqrt{3} & 0 & 0 & 0 & 0 \\
 1 & 0 & 0 & 0 & 0 \\
 -3 & 0 & 0 & 0 & 0 \\
 0 & 0 & 0 & 0 & 0 \\
 1 & 0 & 0 & 0 & 0 \\
\end{array}
\right)
\in\cK.
 $$
So in any realization of this bug $\sB$, if the spheres in $W$ have bends, resp., $a$, $b$, $c$, $d$, and $e$, say, then
$$
0\ = \ (a,b,c,d,e)\cdot K \cdot(1,0,0,0,0)^t 
\ = \ 
-\sqrt{3} a+b-3 c+e.
$$

Suppose this is the case with bends $a,\dots, e$ all integers; then we must have $a=0$.
Note that this applies to not only one particular realization of these five spheres, but to {\it any} such; in particular, the entire orbit of these five spheres under  the symmetry group $\G=\<S\>$ has the bend $a=0$.
 But then $\G$ satisfies an extra polynomial equation, and is not Zariski dense, which is a contradiction.
\epf

\medskip
Lastly, we prove \propref{prop:golden} by exhibiting the following.
\begin{figure}
\includegraphics[width=4in]{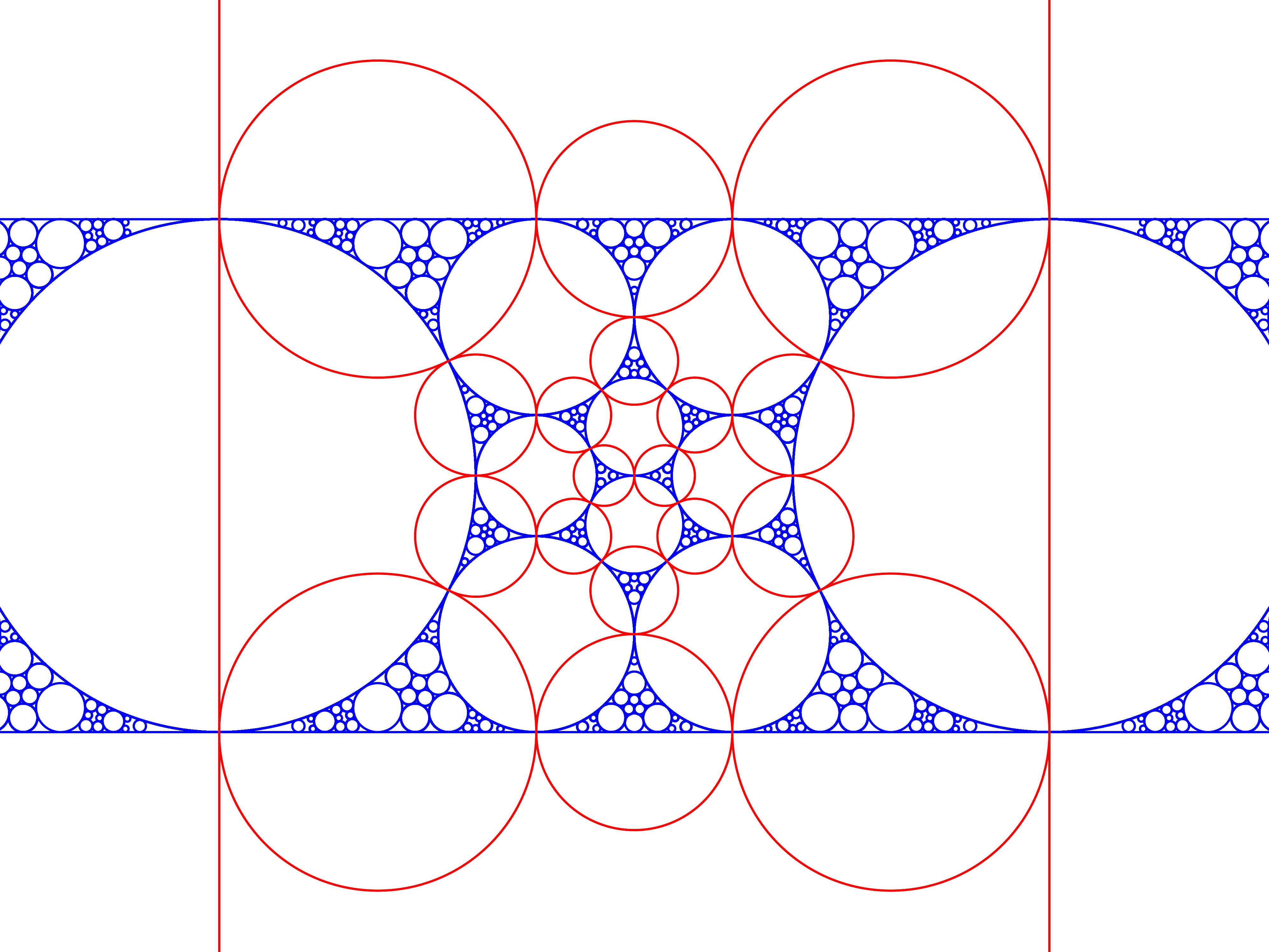}
\caption{Icosahedral packing}
\label{fig:icosa}
\end{figure}
\begin{example}\label{ex:golden}
{\em
In \cite{KN}, a procedure was given for construction a packing $\sP=\sP(\Pi)$ modeled on a (convexly-realizable combinatorial type of a) polyhedron $\Pi$.
In the case of $\Pi$ being the icosahedron, the resulting packing shown in \figref{fig:icosa}. In can be shown that the entire superpacking can be made to have bends in $\fo=\Z[\phi]$, where $\phi={1+\sqrt{5}\over2}$ is the golden mean. But  Vinberg's criterion applied to the ``superGramian'' (in the nomenclature of \cite{KN}) shows that the supergroup is non-arithmetic, because the quadratic form it defines (over $\Q(\phi)$) fails to become definite under the Galois conjugate embedding. 
For details of this computation, see \url{http://math.rutgers.edu/~alexk/maths/Icosahedron.nb}.
}
\end{example}

%{\color{red}Question: McMullen asked about a packing with all spheres disjoint with supergroup commensurate to $\SL_2(\Z[i])$. It seems to me that Millson's construction will give such? Or must the cover  have a  wall tangent to the hyperplane?...}

%\newpage

\section{Hausdorff dimensions of limit sets of packings}\label{sec:Hdim}

%\begin{definition}
%Recall from Definition \ref{def:simple} that a Kleinian bug in $\R^n$ is {\em simple} if it contains no spheres intersecting transversally. 
%\end{definition}

In \propref{prop:Millson} and its application to (the forward direction of) the Classification \thmref{thm:class}, we showed how to construct superintegral Kleinian packings from $\Q$-arithmetic non-uniform lattices. We now extend this construction to 
 show that for infinitely many $j$'s the Hausdorff dimensions of the limit sets of the groups 
$\Ga_j$'s constructed in the proof of \propref{prop:Millson} are pairwise distinct and, moreover, the sequence of Hausdorff dimensions converges to the maximal Hausdorff dimension, $n$.  
This will follow from the following theorem where $\dim$ stands for the Hausdorff dimension, which implies \thmref{thm:abundance}.  We continue with the notation introduced in the proof of \propref{prop:Millson}. 

\begin{thm}\label{thm:monotonic_convergence}
$\lim_{j\to\infty} \dim(\La(\Ga_j))= n$ and for all $j$,  $\dim(\La(\Ga_j))< n$.  
\end{thm}
%\proof 
Note first that the discrete groups $\Ga_j$ contain geometrically finite subgroups which are conjugates of $\pi_1(N)$. In particular,
$$
\delta(\Ga_j)\ge \dim(\La( \pi_1(N)))= n-1\ge \frac{n}{2},
$$
provided that $n\ge 2$.

We next discuss a relation between the Hausdorff dimension of the limit set, the critical exponent and the bottom of the spectrum of the Laplacian. 

For a complete connected Riemannian manifold $M$ let $\la(M)$ denote bottom of the $L^2$-spectrum of the Laplacian of $M$. 
This number can be computed via {\em Rayleigh quotients}: 
$$
\la= \inf \frac{\int_{\Om} |\nabla u|^2}{\int_\Om u^2}
$$
where the infimum is taken over all smooth compactly supported functions (called ``test'' or ``trial'' functions) 
$u\in C^\infty_c(M)$ and $\Om:= \{u>0\}$. 

\begin{thm}
[Semicontinuity of $\la$] $\la(M)$ is upper semicontinuous with respect to the topology of smooth Gromov--Hausdorff convergence:
$$
M_i\to M \Rightarrow \lim\inf_{i\to\infty} \la(M_i) \le \la(M). 
$$
\end{thm}
\proof This is clear using the Rayleigh quotient definition: Every test-function on $M$ is the $C^1$-limit of a sequence of test-functions on $M_i$'s. \qed

\begin{remark}{\em
1. The same theorem (and proof) applies to Riemannian orbifolds. 

2. The bottom of the spectrum is {\em not continuous} with respect to the topology of smooth Gromov--Hausdorff convergence. For instance, for $n\ge 3$ let $M$ be a hyperbolic $n+1$-manifold of finite volume, let $M_i\to M$ be a profinite sequence of (finite) covers of $M$. Then the manifolds $M_i$ converge to the hyperbolic $n$-space, $\la(M_i) =0$, while $\la(\H^{n+1})=n^2/4$. 
}\end{remark}

Given a discrete subgroup $\Ga< \mathrm{Isom}(\H^{n+1})$, let $\delta(\Ga)$ denote the critical exponent of $\Ga$ 
(see e.g. \cite{Nicholls}). A discrete subgroup $\Ga< \mathrm{Isom}(\H^{n+1})$ is called {\em nonelementary} if its limit set consists of more than two points. 
If $\Ga$ is geometrically finite and nonelementary, then $\delta(\Ga)$ equals the Hausdorff dimension of the limit set of $\Ga$, see \cite{Patterson, Sullivan, Tukia, Nicholls}.  We will need the Elstrod--Patterson--Sullivan formula (see e.g. \cite{Nicholls}), relating, for a discrete subgroup $\Ga< \mathrm{Isom}(\H^{n+1})$,  the critical exponent $\delta=\delta(\Ga)$ and the bottom of the spectrum $\la=\la(\H^{n+1}/\Ga)$:

\begin{thm}\label{thm:EPS}
$$
\la=  \left(\frac{n}{2}\right)^2,  \hbox{~~if~~} \delta\le \frac{n}{2}, 
$$
$$
\la= \delta (n- \delta), \hbox{~~if~~} \delta\ge \frac{n}{2}. 
$$
\end{thm}

\medskip

\begin{definition}
A sequence of closed subgroups $\Ga_i$ of a Lie group $G$ is said to converge to a closed subgroup $\Ga< G$ {\em geometrically} or in {\em Chabauty topology} if the following two conditions are met:

1. For every $\ga\in \Ga$ there exists a sequence $\ga_i\in \Ga_i$ which converges to $\ga$. 

2. If a sequence $\ga_{i}\in \Ga_i$ subconverges to $\ga\in G$, then $\ga\in \Ga$.  
\end{definition}

Suppose now that $X$ is a complete connected Riemannian manifold, $G$ is the isometry group of $X$. 
Then $G$ (equipped with the compact-open topology) is a Lie group. Fix a base-point $x\in X$. Consider a sequence of subgroups $\Ga_i< G$ and a subgroup $\Ga< G$ and quotient manifolds/orbifolds $M_i=X/\Ga_i, M=X/\Ga$.  Let $\bar{x}_i, \bar{x}$ denote the projections of $x$ to $M_i$, $M$ respectively.  

The following theorem was proven in \cite{BP} in the case when $X$ is the hyperbolic space  (which will suffice for us), but the same proof works for any complete connected Riemannian manifold. 

\begin{thm}
[See \cite{BP}]
A sequence of discrete subgroups $\Ga_i< G$ converges geometrically to a discrete subgroup $\Ga< G$ if and only if   the sequence of pointed Riemannian manifolds/orbifolds $(M_i, \bar{x}_i)$ 
converges to $(M, \bar{x})$ in the smooth Gromov--Hausdorff topology.   
\end{thm}

\begin{cor}
[Semicontinuity of $\la$ and $\delta$]

(a) Suppose that $\Ga_i$ is a sequence of discrete (nonelementary) subgroups of $\mathrm{Isom}(\H^{n+1})$ converging to  a discrete subgroup $\Ga< \mathrm{Isom}(\H^{n+1})$. Then  
 $$
 \lim\inf_{i\to\infty} \la(\H^{n+1}/\Ga_i) \le \la(\H^{n+1}/\Ga). 
$$

(b) Suppose, in addition, that  $\inf\{\delta(\Ga), \delta(\Ga_i), i\in \N\} \ge n/2$. Then

\begin{equation}\label{eq:lowersemi}
 \lim\inf_{i\to\infty} \delta(\Ga_i) \ge \delta(\Ga). 
\end{equation}
In particular, if $\delta(\Ga_i)\le \delta(\Ga)$ for all $i$ (e.g. if $\Ga_i< \Ga$) then
$$
 \lim_{i\to\infty} \delta(\Ga_i) =\delta(\Ga). 
$$
\end{cor}
\proof The first inequality is a direct corollary of the two theorems above. The second inequality follows 
from the relation of $\delta$ and $\lambda$:  
$$
\la= \delta (n- \delta) \hbox{~~if~~} \delta\ge \frac{n}{2}, 
$$
see \thmref{thm:EPS}. 
\qed

Lastly, we need the following theorem due to T.~Roblin \cite{Roblin} (see also R.~Brooks \cite{Brooks85}): 

\begin{thm}
Let $\Ga< \mathrm{Isom}(\H^{n+1})$ be a geometrically finite\footnote{More generally, 
a group of divergence type, i.e. a group whose Poincar\'e series diverges at the critical exponent.} subgroup and $\hat{\Ga} < \Ga$ a normal subgroup with amenable quotient $\Ga/\hat{\Ga}$.  Then $\delta(\hat\Ga)=\delta(\Ga)$. 
\end{thm}

Thus, if $\Ga$ is a lattice and $\hat{\Ga}< \Ga$ is a normal subgroup with cyclic quotient, then $\delta(\hat \Ga)=n$.  

We can now finish the proof of Theorem \ref{thm:monotonic_convergence}. Let $\Gamma_j$ be the discrete groups as in the theorem. 
 
\begin{lemma}
The sequence $\Ga_j$ geometrically converges to $\hat\Ga$. 
\end{lemma}
\proof This follows from the fact that the sequence of domains $D_j$ exhausts $\hat{M}$ (see  \rmkref{rem:exhaust}), which ensures Gromov--Hausdorff convergence of the corresponding hyperbolic manifolds, hence, geometric convergence of discrete subgroups. \qed 

Thus, we conclude that
$$
\lim_{j\to\infty} \dim(\La(\Ga_j))=\dim(\La(\Ga))=n. 
$$

On the other hand, $\delta(\Ga_j)< n$ since each $\Ga_j$ is geometrically finite and its limit set is a proper subset of $\bS^n$, (Sullivan \cite{Sullivan} and Tukia \cite{Tukia}, independently). This concludes the proof of \thmref{thm:monotonic_convergence} and, hence, of 
 \thmref{thm:abundance}. \qed

\end{document}